\documentclass[12pt]{article}
\usepackage{amsmath,amsfonts,amsthm,mathrsfs}

\newcommand{\R}{{\mathbb R}} \newcommand{\N}{{\mathbb N}}
\newcommand{\K}{{\mathbb K}} \newcommand{\Z}{{\mathbb Z}}
  \def\C{{\mathbb C}}
\newcommand{\Prm}{{\mathbb P}}

\renewcommand{\epsilon}{\varepsilon } 
 
\renewcommand{\rho}{\varrho } 
\renewcommand{\phi}{\varphi }
\newcommand{\E}{{\mathbb E}\,}
\newcommand{\EE}{{\mathbb E}}

\newcommand{\ran}{{\rm ran }}

\newcommand{\ranno}{{\rm ran-non }}
\newcommand{\de}{{\rm det }}
\newcommand{\deno}{{\rm det-non }}
\newcommand{\ca}{{\rm card}}
\newcommand{\avg}{{\rm avg }}
\newcommand{\avgno}{{\rm avg-non }}

\newtheorem{theorem}{Theorem}[section]
\newtheorem{lemma}[theorem]{Lemma}
\newtheorem{corollary}[theorem]{Corollary}
\newtheorem{proposition}[theorem]{Proposition}

\newtheorem{remark}[theorem]{Remark}

\begin {document}
 \title{
Randomized Complexity of Parametric Integration and the Role of Adaption I. Finite Dimensional Case
}

\author {Stefan Heinrich\\
Department of Computer Science\\
RPTU Kaiserslautern-Landau\\
D-67653 Kaiserslautern, Germany}  
\date{\today}
\maketitle

\begin{abstract} 
We study the  randomized $n$-th minimal errors (and hence the complexity) of  vector valued mean computation, which is the discrete version of parametric integration.  The results of the present paper form the basis for the complexity analysis of parametric integration in Sobolev spaces, which will be presented in Part 2. Altogether this extends previous results of Heinrich and Sindambiwe (J.\ Complexity, 15 (1999), 317--341) and Wiegand (Shaker Verlag, 2006). Moreover, a basic problem of Information-Based Complexity on the power of adaption for linear problems in the randomized setting is solved.
\end{abstract}

\section{Introduction}
\label{sec:1}
Let $M,M_1,M_2$ be finite sets and let $1\leq p,q \leq \infty$. We define the space
$L_p^M$ as the set of all functions $f: M \rightarrow \K$
with the norm
\[
\| f \|_{L_p^M} = \left\{\begin{array}{lll}
 \displaystyle \left( \frac{1}{|M|} \sum_{i\in M} |f(i)|^p \right)^{1/p} & \quad\mbox{if}\quad  p<\infty  \\[.5cm]
  \displaystyle\max_{i\in M} |f(i)| & \quad\mbox{if}\quad p=\infty.
    \end{array}
\right.
\]
In the present paper we study the complexity of vector-valued mean computation in the randomized setting.
More precisely, we determine the order of the randomized $n$-th minimal errors of 
\begin{equation}
\label{eq:3}
S^{M_1,M_2} : L_p^{M_1\times M_2} \rightarrow L_q^{M_1}
\end{equation}
with
\begin{equation}
\label{C6}
(S^{M_1,M_2}f)(i)  = \frac{1}{|M_2|} 
\sum_{j\in M_2} f(i,j) .
\end{equation}
The input set is the unit ball of $L_p^{M_1\times M_2}$ and information is standard (values of $f$).
$S^{M_1,M_2}$ can also be viewed as discrete parametric integration. 
For $p=q=\infty$ such an analysis is essentially contained in \cite{HS99} and for $1\le p=q<\infty$ in \cite{Wie06}.

The case $p\ne q$ requires some new techniques. Moreover, it contains a domain of parameters, namely $2<p<q\le \infty$, where adaptive and non-adaptive randomized $n$-th minimal errors deviate by a power of $n$. Since the problem \eqref{C6} is linear, this answers a basic question of Information-Based Complexity (IBC). Let us give some  background on this problem. For a more detailed account on the problems and results around adaption  we refer to \cite{Nov96} and \cite{NW08}, see also \cite{Nov88} and \cite{TWW88}.  

{\bf The adaption problem in the deterministic setting:} It is well-known since the 80ies that for linear problems adaptive and non-adaptive $n$-th minimal errors can deviate at most by a factor of 2, thus for any linear problem $\mathcal{P}=(F,G,S,K,\Lambda)$ (see the definitions below) and any $n\in \N$
\begin{equation}
\label{J9}
e_n^\deno (S,F,G)\le 2e_n^\de (S,F,G),
\end{equation}
see Gal and Micchelli \cite{GM80}, Traub and Wo\'zniakowski \cite{TW80}. Partial results in this direction were shown before by Bakhvalov \cite{Bak71}. Similar results for the average case setting for classes of Gaussian measures were obtained by Wasilkowski and Wo\'zniakowski \cite{WW84}, see also \cite{Was86,Was89}. Kon and Novak \cite{KN90} proved that the  factor 2 in relation \eqref{J9} cannot be replaced by 1. 

{\bf The adaption problem in the randomized setting}: Is there a constant $c>0$ such that for all linear problems
$\mathcal{P}=(F,G,S,K,\Lambda)$ and all $n\in\N$
\begin{equation*}
e_n^\ranno (S,F,G)\le ce_n^\ran (S,F,G) \, {\bf ?}
\end{equation*}
See the open problem on p.\ 213 of \cite{Nov96}, and  Problem 20 on p.\ 146 of \cite{NW08}. Let us note that for some non-linear problems the answer is 'No': for integration of monotone functions \cite{Nov92} and of convex functions \cite{NP94}. (These problems are nonlinear because the input set $F$ is not balanced).

Relations \eqref{A3} and \eqref{A4} of Theorem \ref{theo:1} show: {\bf The answer is 'No' for linear problems.}
The case $2<p<q$ of vector-valued mean computation provides a counterexample.
The paper is organized as follows. In Section \ref{sec:2} we recall the basic notions of IBC and present some auxiliary facts. Moreover, this section contains new general results on the average case setting, which will be needed for the lower bound estimates in the main result. In Section \ref{sec:3} we recall one instant of the randomized norm estimation algorithm from \cite{Hei18} which is a central part of the analysis of the critical domain $2<p<q$. Finally, Section \ref{sec:4} contains the complexity analysis of vector-valued mean computation and the solution of the above mentioned adaption problem.

\section{Preliminaries}
\label{sec:2}

Throughout this paper $\log$ means $\log_2$. We denote $\N=\{1,2,\dots\}$ and $\N_0=\N\cup\{0\}$. The symbol $\K$ stands for the scalar field, which is either $\R$ or $\C$.
We often use the same symbol
$c, c_1,c_2,\dots$ for possibly different constants, even if they appear in a sequence
of relations. However, some constants are supposed to have the same meaning throughout a proof  --  these are denoted by symbols $c(1),c(2),\dots$. The unit ball of a normed space $X$ is denoted by $B_X$.

We work in the framework of IBC \cite{Nov88,TWW88}, using specifically the general approach from \cite{Hei05a, Hei05b}. 
An abstract  numerical problem $\mathcal{P}$ is given as 
\begin{equation}
\label{M7}
\mathcal{P}=(F,G,S,K,\Lambda).
\end{equation}
Here $F$ is a non-empty set, 
$G$ a Banach space and $S$ is a mapping $F\to G$.  The operator $S$ is called the  solution operator, it sends the input  $f\in F$ of our problem to the exact solution $S(f)$. Moreover, $\Lambda$ is a nonempty set of mappings from $F$ to $K$, the set of information functionals, where $K$ is any nonempty set - the set of values of information functionals. 
A problem $\mathcal{P}$ is called linear, if $K=\K$, $F$ is a convex and balanced subset of a linear space $X$ over $\K$,  
$S$ is the restriction to $F$ of a linear operator 
from $X$ to $G$, and each $\lambda\in\Lambda$ is the restriction  to $F$ of a linear mapping from $X$ to $\K$.

A deterministic algorithm for $\mathcal{P}$ 
is a tuple $A=((L_i)_{i=1}^\infty, (\tau_i)_{i=0}^\infty,(\varphi_i)_{i=0}^\infty)$
such that 
$L_1\in\Lambda$, $\tau_0\in\{0,1\}$, $\varphi_0\in G$,
and for $i\in \N$
\begin{equation*}
L_{i+1} : K^i\to \Lambda,\quad
\tau_i:  K^i\to \{0,1\},\quad
\varphi_i:  K^i\to G 
\end{equation*}
are arbitrary mappings, where $K^i$ denotes the $i$-th Cartesian power of $K$.
Given an input $f\in F$, we define $(\lambda_i)_{i=1}^\infty$ with $\lambda_i\in \Lambda$ 
as follows:
\begin{eqnarray*}
\lambda_1=L_1, \quad
\lambda_i=L_i(\lambda_1(f),\dots,\lambda_{i-1}(f))\quad(i\ge 2)\label{RC2}.
\end{eqnarray*}
Define $\ca(A,f)$, the cardinality  of $A$ at input $f$, to be $0$ if $\tau_0=1$. If $\tau_0=0$, let $\ca(A,f)$ be
the first integer $n\ge 1$ with $
\tau_n(\lambda_1(f),\dots,\lambda_n(f))=1 $
if there is such an $n$. If $\tau_0=0$ and no such $n\in \N$ exists, 
put $\ca(A,f)=+\infty$.
We define the output $A(f)$ of algorithm $A$ at input $f$ as
\begin{equation*}
A(f)=\left\{\begin{array}{lll}
\phi_0  & \mbox{if} \quad \ca(A,f)\in \{0,\infty\} \\[.2cm]
\phi_n(\lambda_1(f),\dots,\lambda_n(f))  &\mbox{if} \quad 1\le \ca(A,f)=n<\infty. 
\end{array}
\right.
\end{equation*}
The cardinality of $A$ is defined by
$$
{\rm card}(A,F)=\sup_{f\in F}{\rm card} (A,f)
$$
and the error of $A$ in approximating $S$ by
$$
e(S,A,F,G)=\sup_{f\in F}\|S(f)-A(f)\|_G.
$$
Let $\mathscr{A}^\de(\mathcal{P})$ be the set of all deterministic algorithms for $\mathcal{P}$.
Then the deterministic $n$-th minimal error of $S$ is defined as
\begin{equation*}
e_n^{\rm det } (S,F,G)=\inf_{A\in\mathscr{A}^\de(\mathcal{P}),\, \ca(A,F)\le n  }  
e(S,A,F,G).
\end{equation*}

 A deterministic algorithm is called non-adaptive, if all $L_i$ and $\tau_i$ are constant, in other words, 
\begin{equation}
\label{Q0}
L_i\in \Lambda\quad(i\in \N), \quad \tau_i\in\{0,1\}\quad(i\in \N_0).
\end{equation}
The subset of non-adaptive algorithms in  $\mathscr{A}^\de(\mathcal{P})$ is denoted by $\mathscr{A}^\deno (\mathcal{P})$. Correspondingly, we define the non-adaptive deterministic $n$-th minimal error of $S$ by
\begin{equation*}
e_n^\deno (S,F,G)=\inf_{A\in\mathscr{A}^\deno (\mathcal{P}),\, \ca(A,F)\le n  }  
e(S,A,F,G).
\end{equation*}
Clearly, we always have
\begin{equation*}
e_n^\de (S,F,G)\le e_n^\deno (S,F,G)\quad (n\in\N_0).
\end{equation*}

Below we will consider problems on product structures.
Let $\mathcal{P}$ be an abstract numerical problem  \eqref{M7} and assume that 
\begin{eqnarray*}
F&=&F^{(1)}\times F^{(2)}, \quad K=K^{(1)}\cup K^{(2)},\quad  \Lambda=\Lambda^{(1)}\cup\Lambda^{(2)}, 
\\
F^{(\iota)}&\ne& \emptyset,\quad K^{(\iota)}\ne \emptyset \quad (\iota=1,2), \quad\Lambda^{(1)}\ne \emptyset,\quad \Lambda^{(1)}\cap\Lambda^{(2)}=\emptyset,
\end{eqnarray*}
such that $\Lambda^{(1)}$ consists of mappings into $K^{(1)}$, $\Lambda^{(2)}$ of mappings into $K^{(2)}$, and for all $\lambda\in \Lambda^{(2)}$ we have
$\lambda(f,g)=\lambda(f',g)$ $(f,f'\in F^{(1)},g \in F^{(2)})$,
that is, all $\lambda\in \Lambda^{(2)}$ depend only on $g \in F^{(2)}$ (the $\lambda\in \Lambda^{(1)}$ may depend on both $f$ and $g$). For $\lambda\in \Lambda^{(2)}$ we use both the notation $\lambda(f,g)$ as well as $\lambda(g)$. 

Let $A=((L_i)_{i=1}^\infty, (\tau_i)_{i=0}^\infty,(\varphi_i)_{i=0}^\infty)$ be a deterministic algorithm for $\mathcal{P}$. Given $f\in F^{(1)}$ and $g\in F^{(2)}$, let 
\begin{eqnarray*}
\lambda_1=L_1, \quad
\lambda_i=L_i(\lambda_1(f,g),\dots,\lambda_{i-1}(f,g))\quad(i\ge 2).
\end{eqnarray*}
Define   
\begin{eqnarray*}
\ca_{\Lambda^{(1)}}(A,f,g)&=&|\{k\le \ca(A,f,g):\lambda_k\in \Lambda^{(1)}\}|
\\
\ca_{\Lambda^{(2)}}(A,f,g)&=&|\{k\le \ca(A,f,g):\lambda_k\in \Lambda^{(2)}\}|.
\end{eqnarray*}
Clearly, if $A$ is non-adaptive, these quantities do not depend on $(f,g)$. 
Fix $g\in F^{(2)}$. We define the restricted problem $\mathcal{P}_g=(F^{(1)},G,S_g,K^{(1)},\Lambda_g)$ by setting
\begin{eqnarray}
 S_g: F^{(1)}\to G, \;\; S_g(f)=S(f,g),\quad
\Lambda_{g}=\{\lambda(\,\cdot\,,g):\; \lambda\in \Lambda^{(1)}\}.\label{K6}
\end{eqnarray}
To a  given a deterministic algorithm $A$ for $\mathcal{P}$ and $g\in F^{(2)}$ we will associate an algorithm $A_g$ for the restricted problem $\mathcal{P}_g$. The following result extends Lemma 3 of \cite{Hei17a}  and Proposition 2.1 in \cite{Hei20}.
\begin{lemma}
\label{lem:2}
Let $A$ be a deterministic algorithm for $\mathcal{P}$ and let $g\in F^{(2)}$. Then there is a deterministic algorithm $A_g$ for $\mathcal{P}_g$ such that for all  $f\in F^{(1)}$
\begin{eqnarray}
A_g(f)&=&A(f,g)  \label{I1}
\\
{\rm card}(A_g,f)&=&{\rm card}_{\Lambda^{(1)}}(A,f,g). \label{I0}
\end{eqnarray}
Moreover, if $A$ is non-adaptive, $A_g$ can be chosen to be non-adaptive, as well. In this case \eqref{I0} turns into 
\begin{equation}
{\rm card}(A_g)={\rm card}_{\Lambda^{(1)}}(A). \label{J0}
\end{equation}
\end{lemma}
Except for some minor modifications the proof is the same as that in \cite{Hei20}, we therefore only present the construction of $A_g$ from $A$. \\[.2cm]
{\it Sketch of proof of Lemma \ref{lem:2}}.

Let
$A=((L_i)_{i=1}^\infty, (\tau_i)_{i=0}^\infty,(\varphi_i)_{i=0}^\infty)$ and fix $g\in F^{(2)}$. 
Let $\nu_0\in \Lambda^{(1)}$ be any element.
Given an arbitrary sequence $(y_l)_{l=1}^\infty\in (K^{(1)})^\N$, we define two sequences $(\lambda_i)_{i=1}^\infty\in \Lambda^\N$ and $(z_i)_{i=1}^\infty\in K^\N$ inductively as follows.
Let 
\begin{eqnarray}
\lambda_1&=&L_1\label{I2}\\[.1cm]
z_1&=&  \left\{\begin{array}{lll}
  y_1 & \mbox{if} \quad  \lambda_1\in \Lambda^{(1)}  \\
 \lambda_1(g)&\mbox{if} \quad  \lambda_1\in \Lambda^{(2)}.     
    \end{array}
\right. \notag
\end{eqnarray}
Now let $i\ge 1$, assume that $(\lambda_j)_{j\le i}$ and $(z_j)_{j\le i}$ have been defined, 
let 
\begin{equation*}
l=|\{j\le i:\,\lambda_j\in \Lambda^{(1)}\}|,
\end{equation*}
and set 
\begin{eqnarray}
\lambda_{i+1}&=&L_{i+1}(z_1,\dots,z_i)\label{H6}\\[.1cm]
z_{i+1}&=&  \left\{\begin{array}{lll}
  y_{l+1} & \mbox{if} \quad  \lambda_{i+1}\in \Lambda^{(1)}  \\
 \lambda_{i+1}(g)&\mbox{if} \quad  \lambda_{i+1}\in \Lambda^{(2)}.     
    \end{array}
\right. \notag
\end{eqnarray}
Roughly, this is something like the information $A$ produces, when instead of the values $\lambda(f,g)$ for  $\lambda\in\Lambda^{(1)}$ the consecutive values $y_l$ are inserted.
Let $k_0=0$ and define for $l\in\N$
\begin{equation}
\label{B1}
k_l=\min\{i\in\N:\, i>k_{l-1}, \lambda_i\in \Lambda^{(1)}\},
\end{equation}
($\min\emptyset:=\infty$).

Now we define the functions constituting the algorithm $A_g=(L_{l,g})_{l=1}^\infty, (\tau_{l,g})_{l=0}^\infty,(\varphi_{l,g})_{l=0}^\infty)$ for finite substrings $(y_1,\dots,y_l)$ of the given sequence $(y_l)_{l=1}^\infty$. 
Let $l\in\N_0$ and set
\begin{eqnarray}
L_{l+1,g}(y_1,\dots,y_l)&= & \left\{\begin{array}{lll}
\lambda_{k_{l+1}}(\,\cdot\,,g) & \mbox{if} \quad  k_{l+1}<\infty \\
 \nu_0(\,\cdot\,,g) &\mbox{if} \quad  k_{l+1}=\infty     
    \end{array}
\right. \label{C1}
\\[.1cm]
\tau_{l,g}(y_1,\dots,y_l)&=&\left\{\begin{array}{lll}
0 & \mbox{if} \quad  k_{l+1}<\infty\quad\mbox{and}\quad\tau_i(z_1,\dots,z_i)=0 
\\
&\quad\; \text{ for all } i \text{ with } k_l\le i<k_{l+1}
\\[.1cm]
1 &\mbox{if} \quad  k_{l+1}<\infty\quad\mbox{and}\quad\tau_i(z_1,\dots,z_i)=1 
\\
&\quad\; \text{ for some } i \text{ with } k_l\le i<k_{l+1} 
\\[.1cm]
1 & \mbox{if} \quad  k_{l+1}=\infty
\end{array}
\right. \label{H8}
\\[.1cm]
\phi_{l,g}(y_1,\dots,y_l)&=&\left\{\begin{array}{lll}
\phi_{k_l}(z_1,\dots,z_{k_l}) &\mbox{if} \quad  k_{l+1}<\infty\quad\text{and}\quad\tau_i(z_1,\dots,z_i)=0
\\
&\quad\; \text{ for all } i \text{ with } k_l\le i<k_{l+1} 
\\[.1cm]
\phi_i(z_1,\dots,z_i) & \mbox{if} \quad i \text{ is the smallest idex with } k_l\le i<k_{l+1}
\\
&\quad\; \mbox{and }\tau_i(z_1,\dots,z_i)=1
\\[.1cm]
\phi_0 &\mbox{if} \quad  k_{l+1}=\infty\quad\text{and}\quad\tau_i(z_1,\dots,z_i)=0
\\
&\quad\; \text{ for all } i \text{ with } k_l\le i<\infty .  
    \end{array}
\right. \notag
\end{eqnarray}
Since we defined these functions of finite strings by the help of an infinite string, correctness has to be checked in the sense that for each $l\in\N$ and each sequence $(\tilde{y}_j)_{j=1}^\infty\subset K^{(1)}$ with $y_j=\tilde{y}_j$ for all $j\le l$ the respective values of  
$L_{l+1,g}$, $\tau_{l,g}$, and $\phi_{l,g}$ coincide. But this follows readily from the definitions. 

If $A$ is non-adaptive, then by \eqref{Q0}, 
$L_i\in \Lambda$ and $\tau_i\in\{0,1\}$ for all $i\in \N_0$. Consequently, by \eqref{I2} and \eqref{H6}, $\lambda_i=L_i$, and moreover, by \eqref{B1}, the sequence $(k_l)_{l=0}^\infty$ does not depend on $(y_l)_{l=0}^\infty$. Therefore  \eqref{C1} and \eqref{H8} show that $L_{l+1,g}$ and $\tau_{l,g}$ do not depend on $y_1,\dots,y_l$, thus $L_{l+1,g}\in \Lambda_g$, $\tau_{l,g}\in\{0,1\}$,  hence $A_g$ is non-adaptive, as well, and \eqref{J0} follows.

Finally, the inductive verification of \eqref{I1} and \eqref{I0} is straightforward, but somewhat technical. It follows exactly the line of the respective part of the proof of Proposition 2.1 in \cite{Hei20}.

\qed

A randomized algorithm for $\mathcal{P}$ is a tuple
$
A = (( \Omega, \Sigma, {\mathbb P}),(A_{\omega})_{\omega \in \Omega}),
$
where $(\Omega, \Sigma, {\mathbb P})$ is a probability space and for each $\omega \in \Omega$, $A_{\omega}$ is a deterministic algorithm for $\mathcal{P}$. 
Let $n\in\N_0$. Then $\mathscr{A}^{{\rm ran }}(\mathcal{P})$ stands for the class of 
randomized algorithms $A$ for $\mathcal{P}$ with the following properties: For each $f\in F$ the mapping $\omega\to {\rm card}(A_{\omega},f)$ is $\Sigma$-measurable   
and the mapping $\omega\to A_{\omega}(f)$ is $\Sigma$-to-Borel measurable  and ${\mathbb P}$-almost surely separably valued, i.e.,  there is a separable subspace $G_f$ of $G$ such that 
$\Prm\{\omega:\,A_\omega(f)\in G_f\}=1$.
We define the cardinality of $A\in \mathscr{A}^{{\rm ran }}(\mathcal{P})$ as
$$
\ca(A,F)=\sup_{f\in F}{\mathbb E}\,{\rm card} (A_{\omega},f),
$$
the error as 
$$
e(S,A,F,G)=\sup_{f\in F}{\mathbb E}\,\|S(f)-A_{\omega}(f)\|_G,
$$
and the randomized $n$-th minimal error of $S$ as
\begin{equation*}
e_n^{\rm ran } (S,F,G)=\inf_{A\in\mathscr{A}^{\rm ran }(\mathcal{P})}  
e(S,A,F,G).
\end{equation*}
Considering trivial one-point probability spaces $\Omega=\{\omega\}$ immediately yields 
\begin{equation}
\label{G5}
e_n^{\rm ran } (S,F,G)\le e_n^{\rm det } (S,F,G).
\end{equation}

Similarly to the deterministic case we call a randomized algorithm $(( \Omega, \Sigma, {\mathbb P}), (A_{\omega})_{\omega \in \Omega})$ non-adaptive, if $A_{\omega}$ is non-adaptive for all $\omega \in \Omega$. Furthermore, $\mathscr{A}^\ranno(\mathcal{P})$ is the subset of non-adaptive algorithms in $\mathscr{A}^\ran(\mathcal{P})$.
The non-adaptive randomized $n$-th minimal error of $S$ is given by 
\begin{equation*}
e_n^\ranno (S,F,G)=\inf_{A\in\mathscr{A}^\ranno(\mathcal{P}),\, \ca(A,F)\le n  }  
e(S,A,F,G).
\end{equation*}
Then it holds
\begin{equation}
e_n^\ran (S,F,G)\le e_n^\ranno (S,F,G)\quad (n\in\N_0).\label{B3}
\end{equation}
Moreover, in analogy to \eqref{G5} we have
\begin{equation*}
e_n^\ranno(S,F,G)\le e_n^\deno (S,F,G).
\end{equation*}

We also need the average case setting. For the purposes of this paper we consider it only for measures which are supported by a finite subset of $F$. Then the underlying $\sigma$-algebra is assumed to be $2^F$, therefore no measurability conditions have to be imposed on $S$ and the involved   deterministic algorithms. So let $\mu$ be a  probability measure on $F$ with finite support.
Put
\begin{eqnarray*}
{\rm card} (A, \mu) & = & \int_{F} {\rm card}(A,f )  d \mu(f),
\notag\\
e (S, A, \mu,G) &=& \int_{F} \| S(f) 
	- A(f) \|_G  d \mu(f),
\notag\\
e_n^\avg (S, \mu,G) &=& \inf_{A\in\mathscr{A}^\de(\mathcal{P}): \,{\rm card} (A, \mu)\le n} e (S,A, \mu,G),
\\
e_n^\avgno (S, \mu,G) &=& \inf_{A\in\mathscr{A}^\deno(\mathcal{P}): \,{\rm card} (A, \mu)\le n} e (S,A, \mu,G).
\end{eqnarray*}
Similarly to \eqref{B3} we have 
\begin{equation}
e_n^\avg (S,\mu,G)\le e_n^\avgno (S,\mu,G)\quad (n\in\N_0).\label{B7}
\end{equation}
We use the following well-known results to prove lower bounds. 
\begin{lemma}
\label{Ulem:5}
For every probability measure $\mu$ on $F$ of finite support we have 
\begin{eqnarray*}
e_n^\ran(S,F,G)&\ge& \frac{1}{2}e_{2n}^\avg(S,\mu,G)
\\
e_n^\ranno(S,F,G)&\ge& \frac{1}{2}e_{2n}^\avgno(S,\mu,G).
\end{eqnarray*}
\end{lemma}

Next we prove two general lemmas on the average case. They concern product structures. Let $M\in \N$ and let
for $i=1,\dots, M$, $\mathcal{P}_i=(F_i,G_i,S_i,K_i,\Lambda_i)$ be a numerical problem and $\mu_i$ a probability measure on  $F_i$ whose support is a finite set. We assume that for each $i$ none of the elements of $\Lambda_i$ is constant on $F_i$, that is, 
\begin{equation}
\label{J1}
\text{for all\;} \lambda\in\Lambda_i\; \text{there exist\;}f_1,f_2\in F_i\;\text{with\;}\lambda(f_1)\ne \lambda(f_2).
\end{equation}
Let $1\le q\le \infty$ and let $L_q^M(G_1,\dots, G_M)$ be the space of tuples $(g_i)_{i=1}^M$ with $g_i\in G_i$, endowed with the norm  $\big\|(\|g_i\|)_{i=1}^M\big\|_{L_q^M}$. The coordinate projection of $G$ onto $G_i$ is denoted by $P_i$.
We define the product problem $\mathcal{P}=(F,G,S,K,\Lambda)$ by 
\begin{eqnarray}
&&F=\prod_{i=1}^M F_i,\quad G=L_q^M(G_1,\dots, G_M),\quad K=\bigcup_{i=}^M K_i,
\label{M2}\\
&&S=(S_1,\dots,S_M):F\to G,\quad S(f_1,\dots,f_M)=(S_1(f_1),\dots,S_M(f_M)),\label{M3}
\end{eqnarray}
furthermore, let $\mathscr{F}(F,K)$ denote the set of all mappings from $F$ to $K$, let
\begin{equation}
\label{G9}
\Phi_i:\Lambda_i\to \mathscr{F}(F,K), \quad (\Phi_i(\lambda_i))(f_1,\dots,f_i,\dots, f_M)=\lambda_i(f_i),
\end{equation}
and set
\begin{equation}
\label{N2}
\Lambda=\cup_{i=1}^M\Phi_i(\Lambda_i).
\end{equation}
Note that \eqref{J1} implies
\begin{equation}
\label{J2}
\Phi_i(\Lambda_i)\cap\Phi_j(\Lambda_j)=\emptyset\quad(i\ne j).
\end{equation}
For $1\le i\le M$ we put
\begin{equation}
F'_i=\prod_{1\le j\le M,j\ne i} F_j,\label{J3}
\end{equation}
If $i$ is fixed, we  identify, for convenience of notation, 
\begin{equation}
\label{N3}
F\quad \text{with}\quad F_i\times F_i', \quad f=(f_1,\dots,f_i,\dots,f_M)\in F \quad \text{with}\quad f=(f_i,f'_i),
\end{equation}
 where 
\begin{equation}
f'_i=(f_1,\dots,f_{i-1},f_{i+1},\dots, f_M)\in F_i'.
\label{J4}
\end{equation}

For the following lemma we define
\begin{eqnarray}
\label{UN4}
\mu=\prod_{i=1}^M \mu_i, \quad \mu'_i=\prod_{1\le j\le M,j\ne i} \mu_j.
\end{eqnarray}
\begin{lemma}
\label{Ulem:2}
With the notation above, under assumption \eqref{J1}, we have for each $n\in\N_0$
\begin{equation}
\label{UK3}
e_n^\avg(S,\mu,G)\ge \frac{1}{2}\inf\Bigg\{\big\|(e_{\left\lceil 2n_i\right\rceil}^\avg(S_i,\mu_i,G_i))_{i=1}^M\big\|_{L_q^M}:\,n_i\in\R, n_i\ge 0,\sum_{i=1}^Mn_i\le n\Bigg\}.
\end{equation}
\end{lemma}
\begin{proof}
Let $A=((L_k)_{k=1}^\infty, (\tau_k)_{k=0}^\infty,(\varphi_k)_{k=0}^\infty)$ be a deterministic algorithm for $\mathcal{P}$ with $\ca(A,\mu)\le n$.  
Let $n_i(f)$ be the number of information functionals in $\Phi_i(\Lambda_i)$ called by $A$ at input $f$. Setting $n_i=\EE_\mu n_i(f)$, we have 
\begin{equation}
\label{UK2}
\sum_{i=1}^Mn_i\le n.
\end{equation}
Now we use Lemma \ref{lem:2} for the problem
\begin{equation}
\label{K4}
\mathcal{P}^{(i)}=(F,G_i,P_iS,K,\Lambda)
\end{equation}
and algorithm 
\begin{equation}
\label{K5}
P_iA:=((L_k)_{k=1}^\infty, (\tau_k)_{k=0}^\infty,(P_i\varphi_k)_{k=0}^\infty)
\end{equation}
with 
\begin{eqnarray}
&&F^{(1)}=F_i,\quad F^{(2)}=F_i',\quad K^{(1)}=K_i,\quad  K^{(2)}=\bigcup_{j\ne i} K_j,
\label{K3}\\
&&\Lambda^{(1)}=\Phi_i(\Lambda_i), \quad\Lambda^{(2)}=\bigcup_{j\ne i} \Phi_j(\Lambda_j).\label{J6}
\end{eqnarray}
We conclude that for each $f_i'\in F_i'$ there is a deterministic algorithm $A_{i,f'_i}$ for $\mathcal{P}_{f_i'}^{(i)}$ such that for all $f_i\in F_i$ 
\begin{eqnarray}
A_{i,f'_i}(f_i)&=&P_iA(f_i,f'_i)
\label{N4}\\
\ca(A_{i,f_i'},f_i)&=&\ca_{\Phi_i(\Lambda_i)}(P_iA,f_i,f_i')=\ca_{\Phi_i(\Lambda_i)}(A,f_i,f_i')=n_i(f_i,f_i').
\label{N5}
\end{eqnarray}
Observe that by \eqref{K6}
\begin{equation*}
\mathcal{P}_{f_i'}^{(i)}=(F_i,G_i,(P_iS)_{f_i'},K_i,\Lambda_{f_i'}),
\end{equation*}
moreover, for $f_i\in F_i$ 
\begin{equation*}
(P_iS)_{f_i'}(f_i)=P_iS(f_i,f_i')=S_i(f_i),
\end{equation*}
and, since for  $\lambda_i\in\Lambda_i$ we have $(\Phi_i(\lambda_i))(f_i,f_i')=\lambda_i(f_i)$,
\begin{equation*}
\Lambda_{f_i'}=\{\lambda(\,\cdot\,,f_i'):\; \lambda\in \Phi_i(\Lambda_i)\}
=\{(\Phi_i(\lambda_i))(\,\cdot\,,f_i'):\; \lambda_i\in \Lambda_i\}=\Lambda_i.
\end{equation*}
This implies
\begin{equation}
\label{N6}
\mathcal{P}_{f_i'}^{(i)}=\mathcal{P}_i,
\end{equation}
so $A_{i,f'_i}$ is a deterministic algorithm for $\mathcal{P}_i$.
From \eqref{UN4} and  \eqref{N5} we conclude
\begin{equation}
 \EE_{\mu_i'}\ca(A_{i,f_i'},\mu_i)=\EE_{\mu_i'}\EE_{\mu_i}\ca(A_{i,f_i'},f_i)=\EE_\mu n_i(f_i,f_i') =n_i.\label{K7}
\end{equation}
Now we estimate
\begin{eqnarray}
\label{UK1}
\EE_\mu\|S(f)-A(f)\|&=&\EE_\mu\big\|\big(\|P_iS(f)-P_iA(f)\|_{G_i}\big)_{i=1}^M\big\|_{L_q^M}
\nonumber\\
&\ge&\big\|\big(\EE_\mu\|P_iS(f)-P_iA(f)\|_{G_i}\big)_{i=1}^M\big\|_{L_q^M}.
\end{eqnarray}
Furthermore, \eqref{K7} implies $\mu_i'(\{f_i'\in F_i:\,\ca(A_{i,f_i'},\mu_i)\le 2n_i\})\ge 1/2$, therefore
\begin{eqnarray*}
\lefteqn{\EE_\mu\|P_iS(f)-P_iA(f)\|_{G_i}}
\nonumber\\
&=&\int_{F_i'}\int_{F_i}\|S_i(f_i)-A_{i,f_i'}(f_i)\|_{G_i}d\mu_i(f_i)d\mu_i'(f_i')
=\int_{F_i'}e(S_i,A_{i,f_i'},\mu_i,G_i)d\mu_i'(f_i')
\nonumber\\
&\ge& \int_{\{f_i'\in F_i:\,\ca(A_{i,f_i'},\mu_i)\le 2n_i\}} e(S_i,A_{i,f_i'},\mu_i,G_i)d\mu_i'(f_i')
\ge \frac{1}{2}e_{\left\lceil 2n_i\right\rceil}^\avg(S_i,\mu_i,G_i).
\end{eqnarray*}
Inserting this into \eqref{UK1}, we obtain
\begin{eqnarray*}
\EE_\mu\|S(f)-A(f)\|
&\ge& \frac{1}{2}\big\|\big(e_{\left\lceil 2n_i\right\rceil}^\avg(S_i,\mu_i,G_i)\big)_{i=1}^M\big\|_{L_q^M}.
\end{eqnarray*}
This combined with \eqref{UK2} yields \eqref{UK3}.

\end{proof}
Now consider the case that all $\mathcal{P}_i$ are copies of the same problem $\mathcal{P}_1=(F_1,G_1,S_1,K_1,\Lambda_1)$, and similarly,  $\mu_i=\mu_1$ $(i=1,\dots,M)$.
\begin{corollary}
\label{Ucor:1}
Under these assumptions,
\begin{equation}
\label{UK4}
e_n^\avg(S,\mu,G)\ge 2^{-1-1/q}e_{\left\lceil\frac{4n}{M}\right\rceil}^\avg(S_1,\mu_1,G_1).
\end{equation}
\end{corollary}
\begin{proof}
Let $n_i\in\R$, $n_i\ge 0$ with $\sum_{i=1}^M n_i\le n$. and set $I=\{i:\, n_i\le \frac{2n}{M}\}$, consequently, $|I|\ge \frac{M}{2}$. Hence, for $i\in I$, 
\begin{equation*}
e_{\left\lceil 2n_i\right\rceil}^\avg(S_1,\mu_1,G_1)\ge e_{\left\lceil\frac{4n}{M}\right\rceil}^\avg(S_1,\mu_1,G_1),
\end{equation*}
so Lemma \ref{Ulem:2} gives \eqref{UK4}.

\end{proof}

The next lemma concerns non-adaptive algorithms. We assume the same setting \eqref{M2}--\eqref{J4} as introduced for Lemma \ref{Ulem:2}, except for the definition of $\mu$, which here is given as follows. Let $\nu_i\ge 0$ with $\sum_{i=1}^M \nu_i=1$, let $f_{i,0}'\in F_i'$ be any, but fixed elements, and let
\begin{equation*}
J_i:F_i\to F,\quad J_i(f_i)=(f_i,f_{i,0}')\quad (f_i\in F_i).
\end{equation*}
We define the measure $\mu$ on $F$ by setting for a set $C\subset F$
\begin{equation}
\label{B8}
\mu(C)=\sum_{i=1}^M\nu_i\mu_i(J_i^{-1}(C)),
\end{equation}
thus $\mu$ is a probability measure on $F$ of finite support. 
\begin{lemma}
\label{Ulem:3}
With the notation above and under assumption \eqref{J1} we have for each $n\in\N_0$
\begin{equation}
\label{UN5}
e_n^\avgno(S,\mu,G)\ge M^{-1/q}\min\Bigg\{\sum_{i=1}^M\nu_ie_{n_i}^\avgno(S_i,\mu_i,G_i):\,n_i\in\N_0, n_i\ge 0,\sum_{i=1}^Mn_i\le n\Bigg\}.
\end{equation}
\end{lemma}
\begin{proof}
The proof is similar to that of Lemma \ref{Ulem:2}. Let $A$ be a non-adaptive deterministic algorithm for $\mathcal{P}$ with $\ca(A,F)\le n$.  Let $n_i$ be the number of those non-zero information functionals of $A$ which are from $\Phi_i(\Lambda_i)$. Then 
\begin{equation}
\label{UN8}
\sum_{i=1}^Mn_i\le n.
\end{equation}
We use Lemma \ref{lem:2} again, with the same choice \eqref{K4}--\eqref{J6} and conclusions \eqref{N4}--\eqref{N6}, thus, for each $i$ there is a non-adaptive deterministic algorithm $A_{i,f'_{i,0}}$ for $\mathcal{P}_i$ such that for all $f_i\in F_i$ 
\begin{eqnarray*}
A_{i,f'_{i,0}}(f_i)&=&P_iA(f_i,f'_{i,0})=P_iA(J_i(f_i))
\\
\ca(A_{i,f_{i,0}'})&=&\ca_{\Phi_i(\Lambda_i)}(P_iA)=\ca_{\Phi_i(\Lambda_i)}(A)=n_i.
\end{eqnarray*}
Consequently, using also \eqref{B8},
\begin{eqnarray*}
\int_F\|S(f)-A(f)\|_G d\mu(f)&=&\sum_{i=1}^M \nu_i\int_{F_i}\|S(J_i(f_i))-A(J_i(f_i))\|_Gd\mu_i(f_i)
\nonumber\\
&\ge&M^{-1/q}\sum_{i=1}^M \nu_i\int_{F_i}\|P_iS(J_i(f_i))-P_iA(J_i(f_i))\|_{G_i}d\mu_i(f_i)
\nonumber\\
&=&M^{-1/q}\sum_{i=1}^M \nu_i\int_{F_i}\|S_i(f_i)-A_{i,f'_{i,0}}(f_i)\|_{G_i}d\mu_i(f_i)
\nonumber\\
&\ge &M^{-1/q}\sum_{i=1}^M \nu_ie_{n_i}^\avgno(S_i,\mu_i,G_i),
\end{eqnarray*}
which together with \eqref{UN8} implies \eqref{UN5}.

\end{proof}

Similarly to Corollary \ref{Ucor:1} we obtain for the case $\mathcal{P}_i=\mathcal{P}_1$,   $\mu_i=\mu_1$, $\nu_i=M^{-1}$ $(i=1,\dots,M)$
\begin{corollary}
\label{Ucor:2}
\begin{equation}
\label{U0}
e_n^\avgno(S,\mu,G)\ge 2^{-1}M^{-1/q}e_{\left\lfloor\frac{2n}{M}\right\rfloor}^\avgno(S_1,\mu_1,G_1).
\end{equation}
\end{corollary}
\begin{proof}
Let $n_i\in\N_0$,   $\sum_{i=1}^M n_i\le n$ and define $I=\{i:\, n_i\le \frac{2n}{M}\}$, thus $|I|\ge \frac{M}{2}$. Hence, for $i\in I$, 
\begin{equation*}
e_{n_i}^\avgno(S_1,\mu_1,G_1)\ge e_{\left\lfloor\frac{2n}{M}\right\rfloor}^\avgno(S_1,\mu_1,G_1),
\end{equation*}
so the desired result follows from Lemma \ref{Ulem:3}.

\end{proof}
The types of lower bounds stated in the next lemma are well-known in IBC (see \cite{Nov88,TWW88}). 
For the specific form  presented here we refer, e.g., to \cite{Hei05a},  Lemma 6 for statement (i) and to \cite{Hei18a}, Proposition 3.1 for (ii). Statement (iii) is new and tuned to the topic of this paper. We include the short proof which follows standard lines of IBC.
\begin{lemma}\label{lem:5} 
Assume that  $K=\K$, $F$ is a  subset of a linear space $X$ over $\K$,  
$S$ is the restriction to $F$ of a linear operator 
from $X$ to $G$, and each $\lambda\in\Lambda$ is the restriction to $F$  of a linear mapping from $X$ to $\K$.
Let $\bar{n}\in\N$ and suppose there are $(\psi_i)_{i=1}^{\bar{n}}\subseteq F$
such that the sets $\Lambda_i:=\{\lambda\in \Lambda\,:\, \psi_i(\lambda)\ne 0\}\; (i=1,\dots,\bar{n})$
are mutually disjoint.
Then the following hold for all $n\in\N$ with
$4n<\bar{n}$:

(i) If $\sum_{i=1}^{\bar{n}} \alpha_i \psi_i\in F$ for all sequences $(\alpha_i)_{i=1}^{\bar{n}}\in \{-1,1\}^{\bar{n}}$ and $\mu$ is the distribution of $\sum_{i=1}^{\bar{n}} \epsilon_i \psi_i$, where $\epsilon_i$ are independent Bernoulli random variables
with $\Prm\{\epsilon_i=1\}=\Prm\{\epsilon_i=-1\}=1/2$,  then 
$$
e_n^\avg(S,\mu,G)\ge \frac{1}{2}\min\bigg\{\E\Big\|\sum_{i\in I}\epsilon_iS\psi_i\Big\|_G:\,I\subseteq\{1,\dots,\bar{n}\},\,|I|\ge \bar{n}-2n\bigg\}.
$$

(ii) If $\alpha \psi_i\in F$ for all $1\le i\le\bar{n}$ and $\alpha\in \{-1,1\}$, and $\mu$ is the uniform distribution on the set $\{\alpha \psi_i\,:\,1\le i\le\bar{n},\; \alpha\in \{-1,1\}\}$,
then 
$$
e_n^\avg(S,\mu,G)\ge \frac{1}{2}\min_{1\le i\le \bar{n}} \|S\psi_i\|_G.
$$

(iii) Let $M\in\N$ and let $\big(\bar{I}_k\big)_{k=1}^M$ be disjoint non-empty subsets of $\{1,\dots,\bar{n}\}$ with $|\bar{I}_k|=\bar{n}_k$. Suppose that $\sum_{i\in \bar{I}_k} \alpha_i \psi_i\in F$ for all sequences $(\alpha_i)_{i\in \bar{I}_k}\in \{-1,1\}^{\bar{I}_k}$ and $1\le k\le M$. Then  
\begin{eqnarray*}
  e_n^\deno(S,F,G) 
&\ge& \min\Bigg\{\max\bigg\{\Big\|\sum_{i\in I_k}\alpha_iS\psi_i\Big\|_G:\,(\alpha_i)_{i\in I_k}\in \{-1,1\}^{I_k}, 1\le k\le M\bigg\}:
\\
&&\hspace{2.5cm} n_k\in\N_0, n_k\le \bar{n}_k,\sum_{k=1}^M n_k\le n,I_k\subseteq \bar{I}_k,\,|I_k|= \bar{n}_k-n_k\Bigg\}.
\end{eqnarray*}
\end{lemma}
\noindent{\it Proof of (iii).} 
Let $n\in\N$ with $4n<\bar{n}$ and let $A$ be a deterministic non-adaptive algorithm with $\ca(A,F)\le n$.
Let $\Lambda_A$ be the set of information functionals called by $A$ and let 
$$
I_k=\{i\in \bar{I}_k:\,\Lambda_i \cap \Lambda_A= \emptyset  \},\quad n_k=|\bar{I}_k\setminus I_k|.
$$
Then $n_k\le \bar{n}_k$ and the disjointness of the $\Lambda_i$ implies $\sum_{k=1}^M n_k\le n$. 
Now let 
$$
(\alpha_i)_{i\in \bar{I}_k}\in \{-1,1\}^{\bar{I}_k},\quad 
h_l=\sum_{i\in \bar{I}_k\setminus I_k}\alpha_i \psi_i +(-1)^l \sum_{i\in I_k}\alpha_i \psi_i \quad (l=0,1).
$$
Then $h_0,h_1\in F$, $A(h_0)=A(h_1)$, and by the triangle inequality
\begin{eqnarray*}
\sup_{f\in F}\|S(f)-A(f)\|_G&\ge& \max_{l=0,1}\|S(h_l)-A(h_l)\|_G
\\
&=&\max_{l=0,1}\Bigg\|S\bigg(\sum_{i\in \bar{I}_k\setminus I_k}\alpha_i \psi_i +(-1)^l \sum_{i\in I_k}\alpha_i \psi_i \bigg)-A(h_0)\Bigg\|_G
\ge\Big\|\sum_{i\in I_k}\alpha_iS\psi_i\Big\|_G.
\end{eqnarray*}
Consequently, 
\begin{eqnarray*}
\sup_{f\in F}\|S(f)-A(f)\|_G&\ge& \max\bigg\{\Big\|\sum_{i\in I_k}\alpha_iS\psi_i\Big\|_G:\,(\alpha_i)_{i\in I_k}\in \{-1,1\}^{I_k}, 1\le k\le M\bigg\}
\\
&\ge& \min\bigg\{\max\bigg\{\Big\|\sum_{i\in I_k}\alpha_iS\psi_i\Big\|_G:\,(\alpha_i)_{i\in I_k}\in \{-1,1\}^{I_k}, 1\le k\le M\bigg\}
\\
&&n_k\in\N_0, n_k\le \bar{n}_k,\sum_{k=1}^M n_k\le n,I_k\subseteq \bar{I}_k,\,|I_k|= \bar{n}_k-n_k\bigg\}.
\end{eqnarray*}
It remains to take the infimum over all deterministic non-adaptive algorithms $A$ with \linebreak
$\ca(A,F)\ \le n$.

\qed

We need 
the following well-known procedure of 
``boosting the success probability'', which decreases the failure probability by 
repeating the algorithm a number of times and computing the median of the outputs. The following lemma for $\K=\R$ is essentially contained in  \cite{{Hei01}}, where it was derived in the  setting of quantum computation. We include the short proof for the sake of completeness.

Let  $m\in \N$ and define
$\theta_\K:\K^m\to \K$ as follows. If $\K=\R$, let $\theta_\R$ be the mapping given by the median, that is, if $z_1^*\le  \dots\le z_m^*$ is the non-decreasing rearrangement of $(z_1,\dots,z_m)$, then  
\begin{equation*}
\theta_\R(z_1,\dots,z_m)=\left\{\begin{array}{lll}
   z^*_{(m+1)/2}& \quad\mbox{if}\quad  m \text{ is odd}  \\[.2cm]
  \displaystyle \frac{z_{m/2}^*+z_{m/2+1}^*}{2} & \quad\mbox{if}\quad m \text{ is even}.
    \end{array}
\right. 
\end{equation*}
 If $\K=\C$, then we set 
$$
\theta_\C(z_1,\dots,z_m)=\theta_\R(\Re(z_1),\dots,\Re(z_m))+i \theta_\R(\Im(z_1),\dots,\Im(z_m)).
$$

\begin{lemma}
\label{Ulem:2e}
Let $\zeta_1,\dots,\zeta_m$ be independent, identically distributed $\K$-valued random variables on a probability space $(\Omega,\Sigma,\Prm)$, $z\in \K$, $\epsilon>0$ , and assume that $\Prm\{|z-\zeta_1|_G\le\epsilon\}\ge 3/4$. Then
\begin{equation*}
\Prm\{|z-\theta_\K(\zeta_1,\dots,\zeta_m)|\le c_\K\epsilon\}\ge 1-e^{-m/8},
\end{equation*}
with $c_\R=1$ and $c_\K=\sqrt{2}$.
\end{lemma}
\begin{proof}
Let $\chi_i$ be the indicator function of the set
$
\{|z-\zeta_i|>\epsilon\},
$
thus $\Prm\{\chi_i=1\}\le 1/4$. Hoeffding's inequality, see, e.g., \cite{Pol}, p.\ 191,
yields
\begin{equation}
\label{UN9}
\Prm\left\{\sum_{i=1}^m\chi_i\ge \frac{m}{2}\right\}\le
\Prm\left\{\sum_{i=1}^m(\chi_i-\E\chi_i)\ge \frac{m}{4}\right\}\le e^{-m/8}.
\end{equation}
Define 
\begin{equation*}
\Omega_0=\bigg\{\omega\in\Omega:\,|\{i:\,|z-\zeta_i(\omega)|\le \epsilon\}|> \frac{m}{2} \bigg\},
\end{equation*}
then by \eqref{UN9}, $\Prm(\Omega_0)\ge 1-e^{-m/8}$. Fix $\omega\in\Omega_0$.   
It follows that for $\K=\R$
$$
|z-\theta_\R(\zeta_1(\omega)),\dots,\zeta_m(\omega))|\le \epsilon,
$$
and for $\K=\C$
$$
|\Re(z)-\theta_\R(\Re(\zeta_1(\omega)),\dots,\Re(\zeta_m(\omega)))|\le \epsilon,\quad |\Im(z)-\theta_\R(\Im(\zeta_1(\omega)),\dots,\Im(\zeta_m(\omega)))|\le \epsilon,
$$
and therefore
$$
|z-\theta_\C(\zeta_1(\omega)),\dots,\zeta_m(\omega))|\le \sqrt{2}\epsilon.
$$

\end{proof}

Finally we need some results on Banach space valued random variables.
Given $p$ with $1\leq p \leq 2$, we recall from Ledoux and Talagrand \cite{LT91} that the type $p$ 
constant $\tau_p(X)$ of a Banach space $X$ is the smallest 
$c$ with $0< c \leq +\infty$,
such that for all $n$ and all sequences $(x_i)_{i=1}^n \subset X$,
\begin{equation*}
\E \Big\| \sum_{i=1}^n \varepsilon_i x_i \Big\|^p 
          \leq c^p\sum_{i=1}^n \|x_i \|^p ,
\end{equation*}
where $(\varepsilon_i)$ denotes a sequence of independent symmetric 
Bernoulli random
variables with
$\Prm \{ \varepsilon_i = 1\} =\Prm \{ \varepsilon_i = -1\} =\frac12$.
$X$ is said to be of type $p$ if $\tau_p(X)<\infty$. Trivially, each Banach space is of type 1. Type $p$ implies 
type $p_1$ for all $1\le p_1<p$.
For $1\leq p < \infty$ all $L_p$ spaces  are of type $\min (p,2)$. Moreover, the spaces
$L_p^N$ are of type $\min (p,2)$ uniformly in $N$, that is,
$\tau_{\min (p,2)}(L_p^N) \leq c$. Furthermore, $c_1(\log (N+1))^{1/2}\le \tau_{2}(L_\infty^N) \le c_2(\log (N+1))^{\frac{1}{2}}$.

We will use the following result. The case $p_1=p$ of it is contained in  Proposition 9.11 of \cite{LT91}. The extension  to the case of general $p_1$ is Lemma 2.1 of \cite{Hei08b}.
\begin{lemma}
\label{lem:1} 
Let $1\le p\le 2$, $p\le p_1<\infty$. Then there is a constant $c>0$ such that for each Banach space  $X$  of type $p$, each $n\in \N$
 and each sequence of independent, mean zero $X$-valued random variables
$(\zeta_i)_{i=1}^n$ with $\E\|\zeta_i\|^{p_1}<\infty$ $(1\le i\le n)$ the following holds:
\begin{equation}
\label{A11}
\Bigg(\E \Big\|\sum_{i=1}^n \zeta_i\Big\|^{p_1}\Bigg)^{1/p_1}
\le c\tau_p(X)\Bigg(\sum_{i=1}^n\Big(\E \|\zeta_i\|^{p_1}\Big)^{p/p_1}\Bigg)^{1/p}.
\end{equation}
\end{lemma}
\section{Norm estimation}
\label{sec:3}
A key part of one of the algorithms below will be randomized norm estimation. We use an algorithm from \cite{Hei18}. 
 Let $(Q,\mathcal{Q},\rho)$ be a probability space, let $1\le q< p
\le \infty$, and let $p_1$ be such that
\begin{equation}
\label{RK6}
2<p_1<\infty\quad\text{if}\quad p=\infty \quad\text{and} \quad q=1,
\end{equation}
and
\begin{equation}
\label{RE8a}
\frac{1}{p_1}=1+\frac{1}{p}-\frac{1}{q}\quad\text{if}\quad p<\infty \quad\text{or} \quad q>1.
\end{equation}
For $n\in\N$ define  $A_n^{(1)}=\big(A_{n,\omega}^{(1)}\big)_{\omega\in\Omega}$ by setting 
 for $\omega\in \Omega$ and $f\in L_p(Q,\mathcal{Q},\rho)$ ($=L_p(Q)$ for short)
\begin{eqnarray}
A_{n,\omega}^{(1)}(f)&=&\left(\frac{1}{n}\sum_{i=1}^n|f(\xi_i(\omega))|^q\right)^{1/q},\label{QG4c}
\end{eqnarray}
where $\xi_i$ are independent $Q$-valued random variables on a probability space $(\Omega,\Sigma,\Prm)$ with distribution $\rho$. 

First we recall Lemma 3.1 from \cite{Hei18}.
\begin{lemma} 
\label{Qlem:3}
Let $0<\alpha<\infty$. Then for $x,y\in\R$ with $x,y\ge 0$ and $x+y> 0$  
\begin{equation*}
\min(\alpha,1)\max(x,y)^{\alpha-1}|x-y|\le|x^\alpha-y^\alpha|\le\max(\alpha,1)\max(x,y)^{\alpha-1}|x-y|.
\end{equation*}
Moreover, if $1\le \alpha<\infty$, then
\begin{equation*}
|x-y|\le |x^\alpha-y^\alpha|^{1/\alpha}.
\end{equation*}
\end{lemma}
The following is essentially the upper bound from Proposition 6.3 of \cite{Hei18}. No proof was given there, it was just mentioned there that the (quite technical) proof of Proposition 3.4 simplifies to yield Proposition  6.3. For the sake of completeness we include the full proof here.
\begin{proposition}
\label{Rpro:4}
Let $1\le q<p\le \infty$, and let $p_1$ satisfy \eqref{RK6}--\eqref{RE8a}.  
Then there is 
a constant $c>0$ such that for  all probability spaces $(Q,\mathcal{Q},\rho)$, $f\in L_p(Q)$, and  $n\in\N$ 
\begin{eqnarray*}
\left( \E\left|\|f\|_{L_q(Q)}-A_{n,\omega}^{(1)}(f)\right|^{p_1}\right)^{1/p_1} &\le& cn^{\max\left(1/p-1/q,-1/2\right)}\|f\|_{L_p(Q)}.
\end{eqnarray*}
\end{proposition}
\begin{proof}
Let $u=\min(p_1,2)$. The assumption $q<p$ and \eqref{RK6}--\eqref{RE8a} imply 
\begin{eqnarray}
&&1< u\le2,\quad u\le p_1\le p, 
\label{RK8}\\[.2cm]
&&\frac{1}{u}-1=\max\left(\frac{1}{p_1},\frac{1}{2}\right)-1
=\max\left(\frac{1}{p}-\frac{1}{q},-\frac{1}{2}\right).\label{RK7}
\end{eqnarray}
We have $A_{n,\omega}^{(1)}(af)=|a|A_{n,\omega}^{(1)}(f)$ and $\|af\|_{L_q(Q)}=|a|\|f\|_{L_q(Q)}$ for $a\in\R$, so we can assume w.l.o.g.  $f\in B_{L_p(Q)}$, $f\ne 0$. With the help of Lemma \ref{Qlem:3} we obtain 
\begin{eqnarray*}
|\|f\|_{L_q(Q)}-A_{n,\omega}^{(1)}(f)| 
&=&\left|(\|f\|_{L_q(Q)}^q)^{1/q}-(A_{n,\omega}^{(1)}(f)^q)^{1/q}\right|
\nonumber\\
&\le&\max(\|f\|_{L_q(Q)}^q,A_{n,\omega}^{(1)}(f)^q)^{-(q-1)/q}\big|\|f\|_{L_q(Q)}^q-A_{n,\omega}^{(1)}(f)^q\big|
\nonumber\\
&\le&\|f\|_{L_q(Q)}^{-(q-1)}\ \big|\|f\|_{L_q(Q)}^q-A_{n,\omega}^{(1)}(f)^q\big| \quad(\omega\in\Omega).
\end{eqnarray*}
Consequently, 
\begin{eqnarray}
\left(\E|\|f\|_{L_q(Q)}-A_{n,\omega}^{(1)}(f)|^{p_1}\right)^{1/p_1}&\le&\|f\|_{L_q(Q)}^{-(q-1)}\left(\E\left|\|f\|_{L_q(Q)}^q-\frac{1}{n}\sum_{i=1}^n|f(\xi_i)|^q\right|^{p_1}\right)^{1/p_1}.
\label{RB4}
\end{eqnarray}
Setting
$$
\eta_i=\|f\|_{L_q(Q)}^q-|f(\xi_i)|^q,
$$
we conclude from  \eqref{RK8} and Lemma \ref{lem:1} with $X=\K$ (the scalar field is of type $u$)
\begin{eqnarray*}
&&\left(\E\left|\|f\|_{L_q(Q)}^q-\frac{1}{n}\sum_{i=1}^n|f(\xi_i)|^q\right|^{p_1}\right)^{1/p_1}=\left(\E\bigg|\frac{1}{n}\sum_{i=1}^n \eta_i\bigg|^{p_1}\right)^{1/p_1} 
\notag\\
&\le& cn^{-1}\left(\sum_{i=1}^n \left(\E|\eta_i|^{p_1}\right)^{u/p_1}\right)^{1/u}
=c n^{1/u-1}\left(\E|\eta_1|^{p_1}\right)^{1/p_1}
\le c n^{1/u-1}\left(\E|f(\xi_1)|^{qp_1}\right)^{1/p_1}.
\end{eqnarray*}
Together with \eqref{RB4} we arrive at
\begin{eqnarray}
 \left(\E|\|f\|_{L_q(Q)}-A_{n,\omega}^{(1)}(f)|^{p_1}\right)^{1/p_1} 
&\le&  c n^{1/u-1}\|f\|_{L_q(Q)}^{-(q-1)}\left(\E|f(\xi_1)|^{qp_1}\right)^{1/p_1}.
\label{RB8}
\end{eqnarray}
To go on, we first assume $q=1$. 
Taking into account the second relation of \eqref{RK8} and \eqref{RK7}, inequality   (\ref{RB8}) turns into
\begin{eqnarray*}
\left(\E|\|f\|_{L_q(Q)}-A_{n,\omega}^{(1)}(f)|^{p_1}\right)^{1/p_1} 
&\le&cn^{1/u-1}\left(\E|f(\xi_1)|^{p_1} \right)^{1/p_1}\le c n^{\max\left(1/p-1,-1/2\right)},
\end{eqnarray*}
which concludes the proof for $q=1$.

Now we assume $q>1$, them by \eqref{RE8a}, $p_1<p$. Moreover, defining $v$ by
\begin{equation}
\label{RB7}
\frac{1}{v}+\frac{p_1}{p}=1,
\end{equation}
we have $1\le v<\infty$, and by \eqref{RE8a} and \eqref{RB7}
\begin{eqnarray*}
\frac{1}{p_1v}=\frac{1}{p_1}-\frac{1}{p} =1 -\frac{1}{q},
\end{eqnarray*}
hence 
\begin{eqnarray}
\label{RE7a}
(q-1)p_1v= q.
\end{eqnarray}
Next we show that
\begin{eqnarray}
 \E|f(\xi_1)|^{qp_1}&\le&\left(\E|f(\xi_1)|^{(q-1)p_1v}\right)^{1/v}.
\label{RB0}
\end{eqnarray}
Indeed, if $p<\infty$, \eqref{RB7} and H\"older's inequality give
\begin{eqnarray*}
\E|f(\xi_1)|^{p_1}|f(\xi_1)|^{(q-1)p_1}&\le&\left(\E|f(\xi_1)|^{p}\right)^{p_1/p}
 \left(\E|f(\xi_1)|^{(q-1)p_1v}\right)^{1/v}
\nonumber\\
&\le&\left(\E|f(\xi_1)|^{(q-1)p_1v}\right)^{1/v},
\end{eqnarray*}
while for $p=\infty$ we have  $v=1$ and
\begin{eqnarray*}
\E|f(\xi_1)|^{qp_1}&\le&\|f\|_{L_\infty(Q)}^{p_1}
\E|f(\xi_1)|^{(q-1)p_1}\le \E|f(\xi_1)|^{(q-1)p_1},
\end{eqnarray*}
thus \eqref{RB0} is verified. Furthermore, by \eqref{RE7a},
\begin{eqnarray*}
\left(\E|f(\xi_1)|^{(q-1)p_1v}\right)^{1/v}
=\|f\|_{L_{(q-1)p_1v}(Q)}^{(q-1)p_1}
= \|f\|_{L_q(Q)}^{(q-1)p_1}.
\label{RC0}
\end{eqnarray*}
Together with \eqref{RB0} this implies
\begin{equation*}
\left(\E|f(\xi_1)|^{qp_1}\right)^{1/p_1}\le \|f\|_{L_q(Q)}^{q-1}.
\end{equation*}
Inserting the latter into \eqref{RB8} and using \eqref{RK7} we obtain
\begin{eqnarray*}
\left( \E|\|f\|_{L_q(Q)}-A_{n,\omega}^{(1)}(f)|^{p_1}\right)^{1/p_1} &\le&   cn^{1/u-1}
= cn^{\max\left(1/p-1/q,-1/2\right)}.
\end{eqnarray*}

\end{proof}

\section{Vector Valued Mean Computation}
\label{sec:4}
We refer to the definition of vector valued mean computation $S^{M_1,M_2}$ given in \eqref{eq:3}--\eqref{C6}.
In other words, 
\begin{equation*}
S^{M_1,M_2}f=\frac{1}{|M_2|} \sum_{j\in M_2} f_j
\end{equation*}
is the mean of the vectors 
\begin{equation}
\label{UWN1}
f_j=(f(i,j))_{i\in M_1}.
\end{equation}
It is easily checked by  H\"older's inequality that
\begin{equation}
\label{WN2}
\|S^{M_1,M_2}\|=|M_1|^{\left(1/p-1/q\right)_+},
\end{equation}
(with $a_+:=\max(a,0)$ for $a\in\R$). Expressed in the terminology of Section \ref{sec:2}, see \eqref{M7}, we shall study the problem 
\begin{equation}
\label{AL1}
\left( B_{L_p^{M_1\times M_2}},L_q^{M_1},S^{M_1,M_2},\K,\Lambda\right),
\end{equation}
where
$\Lambda=\{\delta_{ij}:\, i\in M_1,\,j\in M_2\}$ with $\delta_{ij}(f)=f(i,j)$. Clearly, this problem is linear.
For $N_1,N_2\in\N$ we write $L_p^{N_1}$ for $L_p^{\Z[1,N_1]}$, where $\Z[1,N_1]:=\{1,2,\dots,N_1\}$, furthermore
$L_p^{N_1,N_2}$ for $L_p^{\Z[1,N_1]\times \Z[1,N_2]}$, and $S^{N_1,N_2}$ for $S^{\Z[1,N_1],\Z[1,N_2]}$. Along with problem \eqref{AL1}  we will also consider the problem 
\begin{equation}
\label{AL2}
\left( L_p^{M_1\times M_2},L_q^{M_1},S^{M_1,M_2},\K,\Lambda\right),
\end{equation}
i.e., just a larger input space $F$, which is convenient for the  definition and analysis of concrete algorithms . Due to the obvious identifications, it suffices to consider $S^{N_1,N_2}$ for the rest of the paper. 
If $N_1=1$, $S^{N_1,N_2}$ turns into the mean operator $S^{N_2}g=\frac{1}{N_2}\sum_{j=1}^{N_2} g(j)$.

Given  $n\in\N$, $n<N_1N_2$, we define for problem \eqref{AL2} a non-adaptive randomized algorithm 
$$
A_{n}^{(2)}=\big(A_{n,\omega}^{(2)}\big)_{\omega\in\Omega}
$$ 
with
$(\Omega,\Sigma,\mu)$ a suitable probability space as follows. Let $\eta_l \;(l=1,\dots ,\big\lceil\frac{n}{N_1}\big\rceil)$ be independent uniformly distributed on 
$\{ 1, \dots , N_2 \}$ random variables, defined on $(\Omega,\Sigma,\mu)$. 
We put for $f\in L_p^{N_1,N_2}$, $1\le i\le N_1$
\begin{eqnarray}
\big(A_{n,\omega}^{(2)} f\big)(i) &=&0 \quad( n<N_1)
\label{B4}\\
\big(A_{n,\omega}^{(2)} f\big)(i) &=& \left\lceil\frac{n}{N_1}\right\rceil^{-1} \sum_{l=1}^{\left\lceil\frac{n}{N_1}\right\rceil} f(i, \eta_l(\omega))\quad(N_1\le n<N_1N_2).
\label{WM3a}
\end{eqnarray}
\begin{remark}
The constants in the subsequent statements and proofs are independent of the parameters $n$, $N_1$,$N_2$, and $m$. This is also made clear by the order of quantifiers in the respective statements.
\end{remark}

The following result is a slight extension to the case $p\ne q$ of the upper bounds in Wiegand's Theorem 4.2 in \cite{Wie06}.
\begin{proposition}
	\label{pro:4}
	Let $1 \le p,q \le \infty $, $1\le w\le p$, $w<\infty$, and put $\bar{p}=\min(p,2)$. Then there is a constant $c>0$ such that for all $n,N_1,N_2\in \N$ with $n<N_1N_2$ and all $f\in L_p^{N_1,N_2}$ 
	\begin{equation}
	\label{K1}
	\E A_{n,\omega}^{(2)}f= S^{N_1,N_2}f\quad(n\ge N_1),\quad \ca(A_{n,\omega}^{(2)},f)\le 2n \quad(\omega\in\Omega),
	\end{equation}
and
\begin{eqnarray}
	\label{K2}
		\lefteqn{
		 \left(\E \|S^{N_1,N_2}f-A_{n,\omega}^{(2)}f\|_{L_q^{N_1}}^w\right)^{1/w}}
		\nonumber\\[.2cm]
		&\le &
	cN_1^{\left(1/p-1/q\right)_+}\|f\|_{L_p^{N_1,N_2}}	\left\{\begin{array}{lll}
		\left\lceil\frac{n}{N_1}\right\rceil^{-1+1/\bar{p}} & \quad\mbox{if}\quad p<\infty \vee q<\infty   
      \\[.2cm]
		\left\lceil\frac{n}{N_1}\right\rceil^{-1/2} 
		\min\left(\log (N_1+1),\left\lceil\frac{n}{N_1}\right\rceil\right)^{1/2}& \quad\mbox{if}\quad p=q=\infty.
		\end{array}
		\right. 
	\end{eqnarray}
\end{proposition}
\begin{proof}
Relation \eqref{K1} is obvious, while  \eqref{K2} for $n<N_1$  directly follows from \eqref{WN2} and \eqref{B4}. Thus, in the subsequent proof we assume $n\ge N_1$.   Next we prove \eqref{K2} for $p=q$. This case is essentially contained in the proof of Theorem 4.2 in Wiegand \cite{Wie06}, for the sake of completeness we include the short argument. 
For the proof we assume $w=p$ if $p<\infty$ and $2\le w<\infty $ if $p=\infty$. The remaining cases follow from H\"older's inequality.
With $f_j \in L_p^{N_1}$ being defined according to \eqref{UWN1}
we get from \eqref{A11} 
\begin{eqnarray}
    \label{mc-eq:3}
	\lefteqn{\left(\mathbb{E} \|S^{N_1,N_2} f -A_{n,\omega}^{(2)} f \|_{L_p^{N_1}}^{w}\right)^{1/{w}}} \notag \\  
	& = & 
	\left\lceil\frac{n}{N_1}\right\rceil^{-1}\Bigg(\E\Bigg\| \sum_{l=1}^{\lceil n/N_1\rceil} \left( \E f_{\eta_l}-f_{\eta_l} \right)
	\Bigg\|_{L_p^{N_1}}^{w}\Bigg)^{1/{w}}
\notag\\
&\le& c\tau_{\bar{p}}(L_p^{N_1}){\left\lceil\frac{n}{N_1}\right\rceil}^{-1}\Bigg(\sum_{l=1}^{\lceil n/N_1\rceil} \left(\E\Big\|\E f_{\eta_l}-f_{\eta_l}
\Big\|_{L_p^{N_1}}^{w}\right)^{\bar{p}/{w}}\Bigg)^{1/\bar{p}}
\notag\\
&\le& c\tau_{\bar{p}}(L_p^{N_1}){\left\lceil\frac{n}{N_1}\right\rceil}^{-1}\Bigg(\sum_{l=1}^{\lceil n/N_1\rceil} \left(\E\|f_{\eta_l}\|_{L_p^{N_1}}^{w}\right)^{\bar{p}/{w}}\Bigg)^{1/\bar{p}}
\notag\\
&=&c\tau_{\bar{p}}(L_p^{N_1})\left\lceil\frac{n}{N_1}\right\rceil^{-1}\Bigg(\left\lceil\frac{n}{N_1}\right\rceil \left\|\Big(\|f_j\|_{L_p^{N_1}}\Big)_{j=1}^{N_2}\right\|_{L_w^{N_2}}^{\bar{p}}\Bigg)^{1/\bar{p}}\le c\tau_{\bar{p}}(L_p^{N_1}) \left\lceil\frac{n}{N_1}\right\rceil^{-1+1/\bar{p}} \| f \|_{L_p^{N_1,N_2}}
\notag\\
&\le& c\| f \|_{L_p^{N_1,N_2}}\left\{\begin{array}{lll}
\left\lceil\frac{n}{N_1}\right\rceil^{-1+1/\bar{p}} & \quad\mbox{if}\quad 1\le p<\infty    \\[.2cm]
\left\lceil\frac{n}{N_1}\right\rceil^{-1/2} 
\min\left(\log (N_1+1),\left\lceil\frac{n}{N_1}\right\rceil\right)^{1/2}& \quad\mbox{if}\quad p=\infty,
\end{array}
\right. 
\label{C4}
\end{eqnarray}
where the second term in the minimum of \eqref{C4} for $p=\infty$ comes from the  bound
$$
\|S^{N_1,N_2} f -A_{n,\omega}^{(2)} f \|_{L_\infty^{N_1}}\le 2 \| f \|_{L_\infty^{N_1,N_2}}\quad (\omega\in\Omega),
$$
which is an obvious consequence of \eqref{WN2} and  \eqref{WM3a}. This shows  \eqref{K2} for $p=q$.

For $p\ne q$ we have
\begin{equation}
\label{UL4}
\left(\mathbb{E} \|S^{N_1,N_2} f -A_{n,\omega}^{(2)} f \|_{L_q^{N_1}}^{w}\right)^{1/{w}}\le c N_1^{\left(1/p-1/q\right)_+}\left(\mathbb{E} \|S^{N_1,N_2} f -A_{n,\omega}^{(2)} f \|_{L_p^{N_1}}^{w}\right)^{1/{w}}.
\end{equation}
This together with \eqref{C4} gives the desired result except for the case $p=\infty$, $q<\infty$. That case follows by setting $p_1=\max(q,2,w)$ and representing
$$
S^{N_1,N_2}: L_\infty^{N_1,N_2}\overset{J}{\longrightarrow}L_{p_1}^{N_1,N_2}\overset{S^{N_1,N_2}}{\longrightarrow}L_q^{N_1},
$$ 
with $J$ the identical embedding.

\end{proof}
Now we consider the case $2<p<q\le \infty$ and define an adaptive randomized algorithm for problem \eqref{AL2}. Let $f\in L_p^{N_1,N_2}$ and set $f_i=(f(i,j))_{j=1}^{N_2}$ (note that now the $f_i$'s are the rows).  Let $m,n\in\N$  and let
\begin{equation}
\label{L6}
\left\{\xi_{jk}:\, 1\le j\le \left\lceil\frac{n}{N_1}\right\rceil,\,1\le k\le m\right\}, \quad \left\{\eta_{jk}:\, 1\le j\le n,\,1\le k\le m\right\},
\end{equation}
be independent random variables on a probability space $(\Omega,\Sigma,\Prm)$ uniformly distributed over $\{1,\dots,N_2\}$. It is convenient for us to assume that $(\Omega,\Sigma,\Prm)=(\Omega_1,\Sigma_1,\Prm_1)\times (\Omega_2,\Sigma_2,\Prm_2)$, that the $(\xi_{jk})$ are defined on $(\Omega_1,\Sigma_1,\Prm_1)$, and the $(\eta_{jk})$ on $(\Omega_2,\Sigma_2,\Prm_2)$. Furthermore, the expectations with respect to the corresponding  probability spaces are denoted by $\EE,\EE_1,\EE_2$.

We first apply  $k$ times algorithm $A_{\left\lceil n/N_1\right\rceil}^{(1)}$ to estimate $\|f_i\|_{L_2^{N_2}}$ and compute the median of the results . Thus, we put for $\omega_1\in\Omega_1$, $1\le i\le N_1$,  $1\le k\le m$
\begin{eqnarray*}
a_{ik}(\omega_1)&=&\Bigg(\left\lceil\frac{n}{N_1}\right\rceil^{-1}\sum_{1\le j\le \left\lceil\frac{n}{N_1}\right\rceil} |f_i(\xi_{jk}(\omega_1))|^2\Bigg)^{1/2}
\nonumber\\[.2cm]
\tilde{a}_i(\omega_1)&=&\theta_\R\big((a_{ik}(\omega_1))_{k=1}^{m}\big).
\end{eqnarray*}
Next we define the number of samples to be taken in every row, setting for $\omega_1\in\Omega_1$,  $1\le i\le N_1$
\begin{equation}
\label{UM5}
 n_i(\omega_1)= \left\{\begin{array}{lll}
\displaystyle \left\lceil\frac{n}{N_1}\right\rceil& \quad\mbox{if}\quad  \tilde{a}_i(\omega_1)^2\le N_1^{-1}\displaystyle \sum_{l=1}^{N_1}\tilde{a}_l(\omega_1)^2  \\[.5cm]
  \displaystyle\left\lceil\frac{\tilde{a}_i^2n}{\sum_{l=1}^{N_1}\tilde{a}_l^2}\right\rceil & \quad\mbox{if}\quad  \tilde{a}_i(\omega_1)^2> N_1^{-1}\displaystyle \sum_{l=1}^{N_1}\tilde{a}_l(\omega_1)^2 ,   
    \end{array}
\right. 
\end{equation}
and approximate $\big(S^{N_1,N_2}f\big)(i)=S^{N_2}f_i$ for  $\omega_2\in\Omega_2$ by
\begin{eqnarray}
b_{ik}(\omega_1,\omega_2)&=&\frac{1}{n_i(\omega_1)}\sum_{j=1}^{n_i(\omega_1)}f_i(\eta_{jk}(\omega_2)) \quad (1\le k\le m)
\label{UN2}
\\
\tilde{b}_i(\omega_1,\omega_2)&=&\theta_\K \big((b_{ik}(\omega_1,\omega_2))_{k=1}^{m}\big).\label{U1}
\end{eqnarray}
Finally we define the output $A_{n,m,\omega}^{(3)}(f)\in L_q^{N_1} $ of the algorithm for $\omega=(\omega_1,\omega_2)$ as 
\begin{equation}
\label{UN3}
A_{n,m,\omega}^{(3)}(f)=\left\{\begin{array}{lll}
  0 & \quad\mbox{if}\quad n<N_1   \\
 \big(\tilde{b}_i(\omega_1,\omega_2)\big)_{i=1}^{N_1} &  \quad\mbox{if}\quad N_1\le n<N_1N_2.
    \end{array}
\right. 
\end{equation}
\begin{proposition}
\label{pro:5}
Let $2< p <q\le  \infty$ and $1\le w<\infty$. Then there exist constants $c_1,c_2>0$ such that the following hold for all $m,n,N_1,N_2\in\N$ and $ f\in L_p^{N_1,N_2}$:
\begin{equation}
\label{WM1}
\ca(A_{n,m,\omega}^{(3)},f)\le 6mn\quad(\omega\in\Omega)
\end{equation}
and for $m\ge c_1\log(N_1+N_2)$, 
$1 \le n < N_1N_2$ 
\begin{eqnarray}
\label{mc-eq:2}
\left(\E \|S^{N_1,N_2}f-A_{n,m,\omega}^{(3)}f\|_{L_q^{N_1}}^{w}\right)^{1/w}
&\le& 
 c_2\Bigg(N_1^{1/p-1/q}\left\lceil\frac{n}{N_1}\right\rceil^{-\left(1-1/p\right)}+\left\lceil\frac{n}{N_1}\right\rceil^{-1/2}\Bigg)\|f\|_{L_p^{N_1,N_2}}.\quad
\end{eqnarray}
\end{proposition}
\begin{proof}
For $n<N_1$ relation \eqref{WM1} is trivial and \eqref{mc-eq:2} follows from \eqref{WN2}. Hence in the sequel we assume $n\ge N_1$. Note that $1\le n_i\le n$ and the total number of samples is 
\begin{eqnarray*}
mN_1\left\lceil\frac{n}{N_1}\right\rceil+m\sum_{i=1}^{N_1} n_i&\le& 2mN_1\left\lceil\frac{n}{N_1}\right\rceil+
m\sum_{i=1}^{N_1}\left\lceil\frac{\tilde{a}_i^2n}{\sum_{l=1}^{N_1}\tilde{a}_l^2}\right\rceil
\nonumber\\[.2cm]
&\le& 3mn+3mN_1\le 6mn,
\end{eqnarray*}
which is \eqref{WM1}.
Fix $f\in L_p^{N_1,N_2}$. By Proposition \ref{Rpro:4} 
\begin{equation}
\label{C0}
\EE_1\Big|\|f_i\|_{L_2^{N_2}}-a_{ik}\Big|\le c(1) \bigg(\frac{n}{N_1}\bigg)^{-\left(1/2-1/p\right)}\|f_i\|_{L_p^{N_2}}
\end{equation}
 and therefore, 
\begin{equation}
\label{A6}
\Prm_1\left\{\omega_1\in \Omega_1:\,\Big|\|f_i\|_{L_2^{N_2}}-a_{ik}(\omega_1)\Big|\le 4c(1) \bigg(\frac{n}{N_1}\bigg)^{-\left(1/2-1/p\right)}\|f_i\|_{L_p^{N_2}}\right\}\ge \frac{3}{4}.
\end{equation}
Now we set
\begin{equation}
\label{M4}
c(2)=\frac{8(w+1)}{\log e}
\end{equation}
(recall that $\log$ always means $\log_2$), then $m\ge c(2)\log(N_1+N_2)$ implies $e^{-m/8}\le (N_1+N_2)^{-w-1}$.
From \eqref{A6} and Lemma \ref{Ulem:2e} we conclude
\begin{equation}
\label{UL6}
\Prm_1\left\{\omega_1\in \Omega_1:\,\Big|\|f_i\|_{L_2^{N_2}}-\tilde{a}_i(\omega_1)\Big|\le 4c(1) \bigg(\frac{n}{N_1}\bigg)^{-\left(1/2-1/p\right)}\|f_i\|_{L_p^{N_2}}\right\}\ge 1-(N_1+N_2)^{-w-1}.
\end{equation}
Let
\begin{equation}
\label{L0}
\Omega_{1,0}=\left\{\omega_1\in \Omega_1:\,\Big|\|f_i\|_{L_2^{N_2}}-\tilde{a}_i(\omega_1)\Big|\le 4c(1) \bigg(\frac{n}{N_1}\bigg)^{-\left(1/2-1/p\right)}\|f_i\|_{L_p^{N_2}}\quad(1\le i\le N_1)\right\},
\end{equation}
thus
\begin{equation}
\label{L1}
\Prm_1(\Omega_{1,0})\ge 1-(N_1+N_2)^{-w}.
\end{equation}
 Fix $\omega_1\in \Omega_{1,0}$.
Then by \eqref{L0} for all $i$
\begin{equation}
\label{UL7}
\tilde{a}_i(\omega_1)\le c \|f_i\|_{L_p^{N_2}}.
\end{equation}
Consequently,
\begin{equation}
\label{UL9}
\bigg(\frac{1}{N_1}\sum_{i=1}^{N_1}\tilde{a}_i(\omega_1)^2\bigg)^{1/2}\le c\bigg(\frac{1}{N_1}\sum_{i=1}^{N_1}\|f_i\|_{L_p^{N_2}}^2\bigg)^{1/2}
\le c\bigg(\frac{1}{N_1}\sum_{i=1}^{N_1}\|f_i\|_{L_p^{N_2}}^p\bigg)^{1/p}= c\|f\|_{L_p^{N_1,N_2}}.
\end{equation}
It follows from  \eqref{UM5} that
\begin{equation}
\label{UM0}
n_i(\omega_1)\ge \left\lceil\frac{n}{N_1}\right\rceil\quad (1\le i\le N_1).
\end{equation}
Moreover, by \eqref{UN2} and the standard variance estimate for the Monte Carlo method,
\begin{eqnarray}
\EE_2\bigg|\big(S^{N_1,N_2}f\big)(i)-b_{ik}(\omega_1,\omega_2)\bigg|&=&\EE_2\bigg|S^{N_2}f_i-\frac{1}{n_i(\omega_1)}\sum_{j=1}^{n_i(\omega_1)}f_i(\eta_{jk}(\omega_2))\bigg|
\notag\\
&\le& 
n_i(\omega_1)^{-1/2}\|f_i\|_{L_2^{N_2}} \quad (1\le i\le N_1, 1\le k\le m).
\label{UM1}
\end{eqnarray}
Let 
\begin{equation*}
I(\omega_1):=\bigg\{1\le i\le N_1:\, \tilde{a}_i(\omega_1)\le \frac{\|f_i\|_{L_2^{N_2}}}{2}\bigg\},
\end{equation*}
hence from \eqref{L0}
\begin{equation*}
\|f_i\|_{L_2^{N_2}}\le c\bigg(\frac{n}{N_1}\bigg)^{-\left(1/2-1/p\right)}\|f_i\|_{L_p^{N_2}} \quad (i\in I(\omega_1)),
\end{equation*}
which combined with \eqref{UM0} and \eqref{UM1} gives 
\begin{eqnarray}
\EE_2\bigg|\big(S^{N_1,N_2}f\big)(i)-b_{ik}(\omega_1,\omega_2)\bigg|\le 
c\bigg(\frac{n}{N_1}\bigg)^{-\left(1-1/p\right)}\|f_i\|_{L_p^{N_2}}\quad (i\in I(\omega_1),  1\le k\le m) .
\label{UM2}
\end{eqnarray}
Now assume $i\not\in I(\omega_1)$, thus 
\begin{equation}
\label{UM6}
\tilde{a}_i(\omega_1)> \frac{\|f_i\|_{L_2^{N_2}}}{2}.
\end{equation}
We show that 
\begin{eqnarray}
\EE_2\bigg|\big(S^{N_1,N_2}f\big)(i)-b_{ik}(\omega_1,\omega_2)\bigg|\le 
c\bigg(\frac{n}{N_1}\bigg)^{-1/2}\|f\|_{L_p^{N_1,N_2}}\quad (i\not\in I(\omega_1),  1\le k\le m).
\label{UM4}
\end{eqnarray}
Indeed, if $\tilde{a}_i^2\le N_1^{-1}\sum_{l=1}^{N_1}\tilde{a}_l^2$, then by \eqref{UM5},  \eqref{UM1}, \eqref{UM6}, and \eqref{UL9}
\begin{eqnarray*}
\EE_2\bigg|\big(S^{N_1,N_2}f\big)(i)-b_{ik}(\omega_1,\omega_2)\bigg|&\le& 
2\bigg(\frac{n}{N_1}\bigg)^{-1/2}\tilde{a}_i
\le 2\bigg(\frac{n}{N_1}\bigg)^{-1/2}\bigg(\frac{\sum_{l=1}^{N_1}\tilde{a}_l^2}{N_1}\bigg)^{1/2}
\nonumber\\
&\le& c\bigg(\frac{n}{N_1}\bigg)^{-1/2}\|f\|_{L_p^{N_1,N_2}}.
\end{eqnarray*}
On the other hand, if $\tilde{a}_i^2> N_1^{-1} \sum_{l=1}^{N_1}\tilde{a}_l^2$, the same chain of relations yields
\begin{eqnarray*}
\EE_2\bigg|\big(S^{N_1,N_2}f\big)(i)-b_{ik}(\omega_1,\omega_2)\bigg|&\le& 
2\bigg(\frac{\tilde{a}_i^2n}{\sum_{l=1}^{N_1}\tilde{a}_l^2}\bigg)^{-1/2}\tilde{a}_i
=2\bigg(\frac{n}{N_1}\bigg)^{-1/2}\bigg(\frac{\sum_{l=1}^{N_1}\tilde{a}_l^2}{N_1}\bigg)^{1/2}
\nonumber\\
&\le& c\bigg(\frac{n}{N_1}\bigg)^{-1/2}\|f\|_{L_p^{N_1,N_2}}.
\end{eqnarray*}
This proves \eqref{UM4}.

Combining \eqref{UM2} and \eqref{UM4}, we conclude for $1\le i\le N_1$, $1\le k\le m$
\begin{eqnarray*}
\EE_2\bigg|\big(S^{N_1,N_2}f\big)(i)-b_{ik}(\omega_1,\omega_2)\bigg|&\le& 
c(3)\bigg(\frac{n}{N_1}\bigg)^{-\left(1-1/p\right)}\|f_i\|_{L_p^{N_2}}+ c(3)\bigg(\frac{n}{N_1}\bigg)^{-1/2}\|f\|_{L_p^{N_1,N_2}}.
\end{eqnarray*}
Arguing as above \eqref{C0}--\eqref{UL6} and setting $c(4)=4c_\K c(3)$, with $c_\K$ from Lemma \ref{Ulem:2e}, we obtain from \eqref{U1} for $\omega_1\in\Omega_{1,0}$, $m\ge c(2)\log(N_1+N_2)$,  and $1\le i\le N_1$
\begin{eqnarray}
&&\Prm_2\Bigg\{\omega_2\in \Omega_2:\,\Big|\big(S^{N_1,N_2}f\big)(i)-\tilde{b}_i(\omega_1,\omega_2)\Big|
\nonumber\\
&&\quad\le 
c(4)\bigg(\frac{n}{N_1}\bigg)^{-\left(1-1/p\right)}\|f_i\|_{L_p^{N_2}}+ c(4)\bigg(\frac{n}{N_1}\bigg)^{-1/2}\|f\|_{L_p^{N_1,N_2}}\Bigg\}
\ge 1-(N_1+N_2)^{-w-1}.
\label{L2}
\end{eqnarray}
Define 
\begin{eqnarray}
\Omega_{2,0}(\omega_1)&=&\Bigg\{\omega_2\in \Omega_2:\,\Big|\big(S^{N_1,N_2}f\big)(i)-\tilde{b}_i(\omega_1,\omega_2)\Big|
\nonumber\\
&&\le 
c(4)\bigg(\frac{n}{N_1}\bigg)^{-\left(1-1/p\right)}\|f_i\|_{L_p^{N_2}}+ c(4)\bigg(\frac{n}{N_1}\bigg)^{-1/2}\|f\|_{L_p^{N_1,N_2}}\quad (1\le i\le N_1)\Bigg\},
\label{L3}
\end{eqnarray}
thus from \eqref{L2}, for all $\omega_1\in \Omega_{1,0}$
\begin{equation}
\label{N0}
\Prm_2\big(\Omega_{2,0}(\omega_1)\big)\ge 1-(N_1+N_2)^{-w}.
\end{equation}
Now we set
\begin{eqnarray}
\Omega_0&=&\{(\omega_1,\omega_2)\in \Omega: \omega_1\in \Omega_{1,0},\, \omega_2\in\Omega_{2,0}(\omega_1)\}.
\label{L7}
\end{eqnarray}
Since for fixed $f$ all random variables \eqref{L6} take only finitely many values, it follows readily that $\Omega_0\in\Sigma$. Furthermore, from \eqref{L1} and \eqref{N0}, 
\begin{eqnarray}
\Prm(\Omega_0)&=& \int_{\Omega_{1,0}}\Prm_2(\Omega_{2,0}(\omega_1))d\Prm_1(\omega_1)
\nonumber\\
&\ge& (1-(N_1+N_2)^{-w})^2 > 1-2(N_1+N_2)^{-w}.\label{M1}
\end{eqnarray}
It follows from \eqref{L3} and \eqref{L7} that
\begin{eqnarray}
\lefteqn{\big\|S^{N_1,N_2}f-(\tilde{b}_i(\omega))_{i=1}^{N_1}\big\|_{L_q^{N_1}}}
\notag\\
&\le& 
c(4)\bigg(\frac{n}{N_1}\bigg)^{-\left(1-1/p\right)}\Big\|\big(\|f_i\|_{L_p^{N_2}}\big)_{i=1}^{N_1}\Big\|_{L_q^{N_1}}+ c(4)\bigg(\frac{n}{N_1}\bigg)^{-1/2}\|f\|_{L_p^{N_1,N_2}}
\nonumber\\
&\le& c(4)\Bigg(N_1^{1/p-1/q}\bigg(\frac{n}{N_1}\bigg)^{-\left(1-1/p\right)}+ \bigg(\frac{n}{N_1}\bigg)^{-1/2}\Bigg)\|f\|_{L_p^{N_1,N_2}}
\quad(\omega\in\Omega_0).
\label{UN1}
\end{eqnarray}

To estimate the error on $\Omega\setminus \Omega_0$ we note that for all $\omega\in\Omega$
\begin{equation*}
|\tilde{b}_i(\omega)|\le \max_{1\le k\le m}|b_{ik}|\le \max_{1\le j\le N_2}|f(i,j)|\le N_2^{1/p} \|f_i\|_{L_p^{N_2}}
\end{equation*}
and therefore,
\begin{equation}
\big\|(\tilde{b}_i(\omega))_{i=1}^{N_1}\big\|_{L_q^{N_1}}\le N_2^{1/p} \Big\|\big(\|f_i\|_{L_p^{N_2}}\big)_{i=1}^{N_1}\Big\|_{L_q^{N_1}}\le N_1^{1/p-1/q}N_2^{1/p}\|f\|_{L_p^{N_1,N_2}}.
\label{L9}
\end{equation}
Furthermore, by \eqref{WN2},
\begin{equation}
\label{M0}
\big\|S^{N_1,N_2}f\big\|_{L_q^{N_1}}\le N_1^{1/p-1/q} \|f\|_{L_p^{N_1,N_2}}.
\end{equation}
Combining \eqref{M1}, \eqref{L9}, and \eqref{M0}, we conclude
\begin{eqnarray*}
\lefteqn{\left(\int_{\Omega\setminus \Omega_0}\big\|S^{N_1,N_2}f-(\tilde{b}_i(\omega))_{i=1}^{N_1}\big\|_{L_q^{N_1}}^wd\Prm(\omega)\right)^{1/w}}
\\
&\le& N_1^{1/p-1/q}\Big(1+N_2^{1/p}\Big)\Prm(\Omega\setminus \Omega_0)^{1/w}\|f\|_{L_p^{N_1,N_2}}
\\
&\le& 2N_1^{1/p-1/q}\Big(1+N_2^{1/p}\Big)(N_1+N_2)^{-1}\|f\|_{L_p^{N_1,N_2}}
\\
&\le&4N_1^{1/p-1/q}N_2^{1/p}(N_1+N_2)^{-1}\|f\|_{L_p^{N_1,N_2}}\le 4N_1^{1/p-1/q}N_2^{-\left(1-1/p\right)}\|f\|_{L_p^{N_1,N_2}}
\\
&\le&4N_1^{1/p-1/q}\bigg(\frac{n}{N_1}\bigg)^{-\left(1-1/p\right)}\|f\|_{L_p^{N_1,N_2}},
\end{eqnarray*}
the last relation being a consequence of $n< N_1N_2$. Together with \eqref{UN1} this shows \eqref{mc-eq:2}.

\end{proof}
\begin{proposition}
\label{pro:1}
Let $1\le p,q\le \infty$. Then there exist constants $0<c_0<1$, $c_1\dots c_4>0$ such that for each $n,N_1,N_1\in \N$, with $n<c_0N_1N_2$ there exist probability measures $\mu^{(1)},\dots,\mu^{(4)}$ with finite support in $B_{L_p^{N_1,N_2}}$ such that 
\begin{eqnarray}
e_{n}^\avg(S^{N_1,N_2},\mu^{(1)},L_q^{N_1})
&\ge& c_1 \left\lceil\frac{n}{N_1}\right\rceil^{-1/2}\left( \min \left( \log (N_1+1), \left\lceil\frac{n}{N_1}\right\rceil \right) \right)^{\delta_{q,\infty}/2}  ,
\label{UL1}\\
e_{n}^\avg(S^{N_1,N_2},\mu^{(2)},L_q^{N_1})
&\ge&  c_2N_1^{1/p-1/q}\left\lceil\frac{n}{N_1}\right\rceil^{-\left(1-1/p\right)},
\label{UL2}\\
e_{n}^\avg(S^{N_1,N_2},\mu^{(3)},L_q^{N_1})
&\ge& c_3 \left\lceil\frac{n}{N_1}\right\rceil^{-\left(1-1/p\right)},
\label{UL3}\\
e_{n}^{\rm avg-non}\big(S^{N_1,N_2},\mu^{(4)},L_q^{N_1}\big)
&\ge& c_4 N_1^{1/p-1/q}\left\lceil\frac{n}{N_1}\right\rceil^{-1/2}.\label{UQ3}
\end{eqnarray}
\end{proposition}
\begin{proof}
The proofs of \eqref{UL1} and \eqref{UL2} are similar to Wiegand's lower bound proofs of the case $p=q$, see Theorem 4.2 in  \cite{Wie06}. 
For a number $1\le L\le N_2$ we define $L$ disjoint blocks of $\{1,\dots,N_2\}$ by setting
\begin{equation}
\label{UWL7}
D_{j}= \left\{ (j-1)\left\lfloor \frac{N_2}{L}\right\rfloor+1, \dots, 
j \left\lfloor \frac{N_2}{L}\right\rfloor  \right\}  \quad (j=1, \dots, L).
\end{equation}
We have
\begin{equation}
\label{UL0}
  \frac{N_2}{2L}\le  \left\lfloor \frac{N_2}{L}\right\rfloor=|D_{j}|\le \frac{N_2}{L}.
\end{equation}
We set $c_0=\frac1{21}$, let $n\in \N$ be such that 
\begin{equation}
\label{B11}
1\le n< \frac{N_1N_2}{21}
\end{equation}
 and put 
\begin{equation}
\label{UK5}
L=\left\lfloor\frac{4n}{N_1}\right\rfloor +1,
\end{equation}
hence 
\begin{equation}
\label{B6}
\left\lceil\frac{4n}{N_1}\right\rceil< L\le 5\left\lceil\frac{n}{N_1}\right\rceil
\end{equation}
and, since by \eqref{B11}, $\frac{4n}{N_1}<N_2$, 
\begin{equation*}
 L\le N_2.
\end{equation*}

To prove \eqref{UL1}, we define functions $\psi_{ij}\in L_p^{N_1,N_2}$ by
\begin{equation*}
\psi_{ij}(s,t) = \left\{ \begin{array}{ll}
                1, & \text{ if} \enspace s=i\quad\mbox{and}\quad t\in D_{j}\\
                0  & \text{ otherwise.}
                         \end{array}
                 \right.
\end{equation*}
By the construction of the $\psi_{ij}$,
\[
\sum_{i=1}^{N_1} \sum_{j=1}^L \alpha_{ij} \psi_{ij}\in  B_{L_p^{N_1,N_2}}
\]
for all $\alpha_{ij} = \pm 1$.  Let $(\varepsilon_{ij})_{i=1,j=1}^{N_1,L}$ be independent 
symmetric Bernoulli random variables and let $\mu^{(1)}$ be the distribution of 
$
\sum_{i=1,j=1}^{N_1,L}\varepsilon_{ij}\psi_{ij}.
$
Since by (\ref{B6}), $LN_1>4n$, we can apply Lemma \ref{lem:5}. So let $\mathcal{K}$ be any subset of 
$\{(i,j)\,:\,1\le i\le N_1,\, 1\le j\le L\}$ with $|\mathcal{K}|\ge LN_1-2n$. Then 
\begin{equation*}
|\mathcal{K}|\ge\frac12 LN_1.
\end{equation*}
For $1\le i\le N_1$ let 
$$
\mathcal{K}_i=\{1\le j\le L\,:\, (i,j)\in \mathcal{K}\}
$$
and
\[
I := \left\{ 1\le i\le N_1 \, :\, |\mathcal{K}_i| \geq \frac{L}{4} \right\} .
\]
Then 
\begin{equation}
\label{UWL4}
|I| \geq \frac{N_1}{4} .
\end{equation}
Let $(e_i)_{i=1}^{N_1}$ denote the unit vectors in  $\K^{N_1}$,  $(g_i)_{i=1}^{\left\lceil N_1/4\right\rceil}$ the unit vectors in  $\K^{\left\lceil N_1/4\right\rceil}$.
Then from \eqref{UL0}, \eqref{UWL4}, and the contraction principle for Rademacher series
(see \cite{LT91}, Theorem 4.4) we get
\begin{eqnarray*}
\lefteqn{\E \bigg\|   \sum_{(i,j)\in \mathcal{K}} \varepsilon_{ij} S^{N_1,N_2}\psi_{ij}  \bigg\|_{L_q^{N_1}}}
\\
 &\ge& \frac{|D_1|}{N_2}\E \bigg\|  \sum_{i\in I} \sum_{j\in \mathcal{K}_i} \varepsilon_{ij}e_{i}   \bigg\|_{L_q^{N_1}}=\frac{|D_1|}{N_2}\E \bigg\|  \sum_{i=1}^{|I|} \sum_{j=1} ^{|\mathcal{K}_i|} \varepsilon_{ij}e_{i}   \bigg\|_{L_q^{N_1}}
\nonumber\\
&\ge& \frac{|D_1|}{N_2}\E \Bigg\| \sum_{i=1}^{\left\lceil N_1/4\right\rceil}\sum_{j=1}^{\left\lceil L/4\right\rceil} \varepsilon_{ij}e_{i}   \Bigg\|_{L_q^{N_1}} = \frac{|D_1|\left\lceil N_1/4\right\rceil^{1/q}}{N_2N_1^{1/q}}\E \Bigg\| \sum_{i=1}^{\left\lceil N_1/4\right\rceil}\sum_{j=1}^{\left\lceil L/4\right\rceil} \varepsilon_{ij}g_{i}   \Bigg\|_{L_q^{\left\lceil N_1/4\right\rceil}}
\nonumber\\
&\ge&\frac{1}{8L}\E \Bigg\| \sum_{i=1}^{\left\lceil N_1/4\right\rceil}\sum_{j=1}^{\left\lceil L/4\right\rceil} \varepsilon_{ij}g_{i}   \Bigg\|_{L_q^{\left\lceil N_1/4\right\rceil}}
\end{eqnarray*}
and from Lemma \ref{lem:5} (i)
\begin{eqnarray}
\lefteqn{e_n^\avg(S^{N_1,N_2},\mu^{(1)},L_q^{N_1}) }
\nonumber\\
&\ge& \frac12\min_{|\mathcal{K}|\ge LN_1-2n}\E \Bigg\|   \sum_{(i,j)\in \mathcal{K}} \varepsilon_{ij} S^{N_1,N_2}\psi_{ij}  \Bigg\|_{L_q^{N_1}}
\ge  \frac{1}{16L}\E \Bigg\| \sum_{i=1}^{\left\lceil N_1/4\right\rceil}\sum_{j=1}^{\left\lceil L/4\right\rceil} \varepsilon_{ij}g_{i}   \Bigg\|_{L_q^{\left\lceil N_1/4\right\rceil}}.\quad\label{B2}
\end{eqnarray}
For $q=\infty$ we use Lemma 5.3 of \cite{HS99} and \eqref{B6} 
to get
\begin{eqnarray*}
e_n^\avg(S^{N_1,N_2},\mu^{(1)},L_\infty^{N_1})
&\ge&\frac{c}{L}\left(  \left\lceil \frac{L}{4}\right\rceil \min \left(\log \left(\left\lceil \frac{N_1}{4}\right\rceil +1\right),\left\lceil \frac{L}{4}\right\rceil \right) \right)^{1/2}
\notag\\
& \ge & c \left\lceil\frac{n}{N_1}\right\rceil^{-1/2}\min\left(\log (N_1+1),\left\lceil\frac{n}{N_1}\right\rceil\right)^{1/2} . 
\end{eqnarray*}
If $1\le q<\infty$, we denote $\bar{\varepsilon}_{j}=(\varepsilon_{ij})_{i=1}^{\left\lceil N_1/4\right\rceil}\in L_q^{\left\lceil N_1/4\right\rceil}$ and let
$(\alpha_j)_{j=1}^{\left\lceil L/4\right\rceil}$ be independent, also of $\epsilon_{ij}$, symmetric Bernoulli random variables. Then, 
using the equivalence of
moments for Rademacher series and Khintchine's inequality (see \cite{LT91}, Theorem 4.7 and Lemma 4.1)
we get from \eqref{B6}  and \eqref{B2},
\begin{eqnarray*}
\lefteqn{e_n^\avg(S^{N_1,N_2},\mu^{(1)},L_q^{N_1})}
\notag\\
&\ge& \frac{1}{16L}\E \Bigg\| \sum_{j=1}^{\left\lceil L/4\right\rceil} \bar{\varepsilon}_{j}   \Bigg\|_{L_q^{\left\lceil N_1/4\right\rceil}}
=\frac{1}{16L}\E^{(\alpha)}\E^{(\epsilon)} \Bigg\| \sum_{j=1}^{\left\lceil L/4\right\rceil} \alpha_j\bar{\varepsilon}_{j}   \Bigg\|_{L_q^{\left\lceil N_1/4\right\rceil}}
\notag\\[.2cm]
&=&\frac{1}{16L}\E^{(\epsilon)}\E^{(\alpha)} \Bigg\| \sum_{j=1}^{\left\lceil L/4\right\rceil} \alpha_j\bar{\varepsilon}_{j}   \Bigg\|_{L_q^{\left\lceil N_1/4\right\rceil}}
\ge \frac{c}{L}\E^{(\epsilon)}\Bigg( \E^{(\alpha)}\Bigg\| \sum_{j=1}^{\left\lceil L/4\right\rceil} \alpha_j\bar{\varepsilon}_{j}   \Bigg\|_{L_q^{\left\lceil N_1/4\right\rceil}}^q\Bigg)^{1/q}
\notag\\[.2cm]
&=&\frac{c}{L}\E^{(\epsilon)}\Bigg( \left\lceil  \frac{N_1}{4}\right\rceil^{-1}\sum_{i=1}^{\left\lceil N_1/4\right\rceil}\E^{(\alpha)}\Bigg| \sum_{j=1}^{\left\lceil L/4\right\rceil} \alpha_j\epsilon_{ij}   \Bigg|^q\Bigg)^{1/q}\ge\frac{c}{L}\left\lceil\frac{ L}{4}\right\rceil^{1/2}\ge c L^{-1/2} \ge c \left\lceil\frac{n}{N_1}\right\rceil^{-1/2}.
\end{eqnarray*}
This proves \eqref{UL1}.

To show the second lower bound, \eqref{UL2}, we use the same set of blocks 
$D_{j} \;( j=1, \dots, L)$  as defined in (\ref{UWL7}) and the same $L$ given by \eqref{UK5}, put
\begin{equation*}
\psi_{ij}(s,t) = \left\{ \begin{array}{ll}
                N_1^{1/p}N_2^{1/p}|D_{j}|^{-1/p} & \text{ if} \enspace s=i\quad\mbox{and}\quad t\in D_{j},\\
                0  & \text{ otherwise,}
                         \end{array}
                 \right.
\end{equation*}
and let $\mu^{(2)}$ be the uniform distribution on the set 
$$
\{ \alpha \psi_{ij}:\,i=1, \dots, N_1,\, j=1, \dots, L ,\, \alpha =\pm 1\} \subset B_{L_p^{N_1,N_2}} .
$$
Recall that by (\ref{B6}), $LN_1>4n$, so from Lemma \ref{lem:5}(ii) and relations \eqref{UL0} and \eqref{B6} we conclude
\begin{eqnarray*}
e_n^\avg(S^{N_1,N_2},\mu^{(2)},L_q^{N_1})
&\ge& \frac12\big\| S^{N_1,N_2} \psi_{1,1}\big\|_{L_q^{N_1}}
 =  \frac12 N_1^{1/p-1/q}N_2^{-\left(1-1/p\right)}|D_{j}|^{1-1/p} 
\\&\geq& \frac12 N_1^{1/p-1/q}(2L)^{-\left(1-1/p\right)}
\ge cN_1^{1/p-1/q}\left\lceil\frac{n}{N_1}\right\rceil^{-\left(1-1/p\right)},
\end{eqnarray*}
thus \eqref{UL2}.

For the proof of the remaining inequalities \eqref{UL3} and \eqref{UQ3} we can assume $n\ge N_1$, because for $n<N_1$ the already shown relation \eqref{UL1} implies \eqref{UL3}, while \eqref{UL2} together with \eqref{B7} gives \eqref{UQ3}. We set 
\begin{equation}
\label{UV1}
L=4\left\lceil\frac{4n}{N_1}\right\rceil+1,
\end{equation}
hence by \eqref{B11}
$$
L\le \frac{16n}{N_1}+5\le \frac{21n}{N_1}\le N_2.
$$

To prove \eqref{UL3}, 
we apply Corollary \ref{Ucor:1}, where we put
\begin{eqnarray}
 M=N_1, \quad F_1=L_p^{N_2},\quad S_1=S^{N_2},\quad G_1=K_1=\K, \quad\Lambda_1=\{\delta_j:\, 1\le j\le N_2\}
\label{L5}
\end{eqnarray}
with $\delta_j(g)=g(j)$. Then obviously \eqref{J1} is satisfied and
\begin{eqnarray}
&&F=\prod_{i=1}^{N_1} L_p^{N_2}=  L_p^{N_1,N_2},\quad G=L_q^{N_1},\quad S=S^{N_1,N_2},
\label{L4}\\
&& K=\K,\quad \Lambda=\{\delta_{ij}:\, 1\le i\le N_1,\, 1\le j\le N_2\}. 
\label{L8}
\end{eqnarray}
Again we use the blocks 
$D_{j} \;( j=1, \dots, L)$  given by  (\ref{UWL7}) and define $\psi_j\in B_{L_p^{N_2}}$ by 
\begin{equation*}
\psi_{j}(t) = \left\{ \begin{array}{ll}
                N_2^{1/p}|D_{j}|^{-1/p} & \text{ if} \enspace \quad t\in D_{j},\\
                0  & \text{ otherwise.}
                         \end{array}
                 \right.
\end{equation*}
Let $\mu_1$ be the uniform distribution on $\{\alpha\psi_j:\, 1\le j\le L,\, \alpha=\pm 1 \}$. The measure $\mu^{(3)}=\mu_1^{N_1}$, compare \eqref{UN4}, has its support in $B_{L_p^{N_1,N_2}}$ and we derive from Corollary \ref{Ucor:1} 
\begin{equation}
\label{B9}
e_n^\avg(S^{N_1,N_2},\mu^{(3)},L_q^{N_1})\ge 2^{-1-1/q}e_{\left\lceil\frac{4n}{N_1}\right\rceil}^\avg(S^{N_2},\mu_1,\K).
\end{equation}
By Lemma \ref{lem:5}(ii), \eqref{UL0}, and \eqref{UV1}
\begin{equation*}
e_{\left\lceil\frac{4n}{N_1}\right\rceil}^\avg(S^{N_2},\mu_1,\K)\ge \frac{1}{2}|S^{N_2}\psi_1|=N_2^{1/p-1}|D_{j}|^{1-1/p}\ge N_2^{1/p-1} \Bigg(\frac{N_2}{2L}\Bigg)^{1-1/p}\ge c \left\lceil\frac{n}{N_1}\right\rceil^{-\left(1-1/p\right)},
\end{equation*}
which together with \eqref{B9} gives \eqref{UL3}.

Finally, we turn to \eqref{UQ3}, where we use Corollary \ref{Ucor:2} with the same choice \eqref{L5}. Consequently, \eqref{J1}, \eqref{L4}, and \eqref{L8} hold. We set 
\begin{equation}
\label{UP3}
\psi_j=N_1^{1/p}\chi_{D_j}\in L_p^{N_2} \quad (j=1,\dots,L),
\end{equation}
with $D_{j}$ given by  (\ref{UWL7}) and $L$ by  (\ref{UV1}). 
Let $(\epsilon_j)_{j=1}^L$ be independent symmetric Bernoulli random variables, let $\mu_1$ be the distribution of  $ \sum_{j=1}^L  \epsilon_j\psi_j$, and $f_{i,0}'=0$ $(i=1,\dots,N_1)$. Denote the resulting from \eqref{B8} measure by $\mu^{(4)}$. Observe that by \eqref{UP3} $\mu^{(4)}$ is supported by $B_{L_p^{N_1,N_2}}$. 
Now  \eqref{U0} and \eqref{B7} yield
\begin{equation}
\label{UP2}
e_n^\avgno(S^{N_1,N_2},\mu^{(4)},L_q^{N_1})\ge \frac{1}{2} N_1^{-1/q}e_{\left\lceil\frac{2n}{N_1}\right\rceil}^\avgno(S^{N_2},\mu_1,\K)\ge \frac{1}{2} N_1^{-1/q}e_{\left\lceil\frac{2n}{N_1}\right\rceil}^\avg(S^{N_2},\mu_1,\K).
\end{equation}
By Lemma \ref{lem:5}(i), \eqref{UL0}, \eqref{UV1}, \eqref{UP3}, and Khintchine's inequality
\begin{eqnarray*}
e_{\left\lceil\frac{2n}{N_1}\right\rceil}^\avg(S^{N_2},\mu_1,\K)&\ge& \frac{1}{2}\min\bigg\{\E\Big|\sum_{i\in I}\epsilon_iS^{N_2}\psi_i\Big|:\,I\subseteq\{1,\dots,L\},\,|I|\ge L-2\left\lceil\frac{2n}{N_1}\right\rceil\bigg\}
\\
&\ge& c L^{1/2}|S^{N_2}\psi_1|\ge cN_1^{1/p}L^{-1/2}\ge cN_1^{1/p}\left\lceil\frac{n}{N_1}\right\rceil^{-1/2}.
\end{eqnarray*}
Inserting this into \eqref{UP2} finally yields \eqref{UQ3}.

\end{proof}

\begin{theorem}
\label{theo:1}
Let $1 \leq p,q \le \infty $ and put $\bar{p}=\min(p,2)$. Then there exists constants $0<c_0<1$, $c_1,\dots,c_7>0$,  such that for $n,N_1,N_2\in\N$  with 
$n < c_0N_1N_2$ the following hold: \\
If $p\le 2$ or $p\ge q$, then  
\begin{eqnarray}
\label{A2}
&&
c_1N_1^{\left(1/p-1/q\right)_+}\left\lceil\frac{n}{N_1}\right\rceil^{-\left(1-1/\bar{p}\right)}\left( \min \left( \log (N_1+1), \left\lceil\frac{n}{N_1}\right\rceil \right) \right)^{\delta_{p,\infty}\delta_{q,\infty}/2}
\nonumber\\[.2cm]
&\le& 
e_n^\ran(S^{N_1,N_2},B_{L_p^{N_1,N_2}},L_q^{N_1})\le e_n^{\rm ran-non}(S^{N_1,N_2},B_{L_p^{N_1,N_2}},L_q^{N_1}) \nonumber\\[.2cm]
&\le &
c_2N_1^{\left(1/p-1/q\right)_+}\left\lceil\frac{n}{N_1}\right\rceil^{-\left(1-1/\bar{p}\right)}\left( \min \left( \log (N_1+1), \left\lceil\frac{n}{N_1}\right\rceil \right) \right)^{\delta_{p,\infty}\delta_{q,\infty}/2}.
\end{eqnarray}
If $2<p<q$, then  
\begin{eqnarray}
\label{A3}
&&c_3N_1^{1/p-1/q}\left\lceil\frac{n}{N_1}\right\rceil^{-\left(1-1/p\right)}+c_3\left\lceil\frac{n}{N_1}\right\rceil^{-1/2}  \left( \log (N_1+1) \right)^{\delta_{q,\infty}/2}
\nonumber\\[.2cm]
&\le& e_n^\ran(S^{N_1,N_2},B_{L_p^{N_1,N_2}},L_q^{N_1}) 
\nonumber\\[.2cm]
&\le &
 c_4N_1^{1/p-1/q}\left\lceil\frac{n}{N_1\log(N_1+N_2)}\right\rceil^{-\left(1-1/p\right)}+c_4\left\lceil\frac{n}{N_1\log(N_1+N_2)}\right\rceil^{-1/2}
\end{eqnarray}
and
\begin{eqnarray}
c_5N_1^{1/p-1/q}\left\lceil\frac{n}{N_1}\right\rceil^{-1/2}
  &\le& e_n^{\rm ran-non}(S^{N_1,N_2},B_{L_p^{N_1,N_2}},L_q^{N_1})
\le c_6N_1^{1/p-1/q}\left\lceil\frac{n}{N_1}\right\rceil^{-1/2}.
\label{A4}
\end{eqnarray}
In the deterministic setting we have 
\begin{eqnarray}
\label{AC1}
c_7N_1^{\left(1/p-1/q\right)_+}
&\le& 
e_n^\de(S^{N_1,N_2},B_{L_p^{N_1,N_2}},L_q^{N_1})
\nonumber\\
&\le& e_n^\deno(S^{N_1,N_2},B_{L_p^{N_1,N_2}},L_q^{N_1})
\le 
N_1^{\left(1/p-1/q\right)_+}.
\end{eqnarray}

\end{theorem}
\begin{proof}
First we mention that for all lower bounds we use the relation between average case and randomized setting, Lemma \ref{Ulem:5}, without further notice.
For $1\le n<N_1$ the upper bounds follow from \eqref{WN2}, the lower bounds from \eqref{UL2} and \eqref{UL3} of Proposition \ref{pro:1}. 
In the sequel we assume $n\ge N_1$. The upper bounds in \eqref{A2} and \eqref{A4} are a consequence of Proposition \ref{pro:4}, since the involved algorithm is non-adaptive. If $n<6N_1\left\lceil c(1)\log(N_1+N_2)\right\rceil$, where $c(1)$ stands for the constant $c_1$ from Proposition \ref{pro:5}, the upper bound of \eqref{A3} follows from \eqref{WN2}. Now assume
\begin{equation}
\label{C11}
n\ge 6N_1\left\lceil c(1)\log(N_1+N_2)\right\rceil.
\end{equation}
We set 
\begin{equation*}
m=\left\lceil c(1)\log(N_1+N_2)\right\rceil, \quad\tilde{n}=\left\lfloor\frac{n}{6\left\lceil c(1)\log(N_1+N_2)\right\rceil}\right\rfloor,
\end{equation*}
and use Proposition \ref{pro:5} with $\tilde{n} $ instead of $n$. Hence by \eqref{WM1}
\begin{equation*}
\ca(A_{\tilde{n},m,\omega}^{(3)},L_p^{N_1,N_2})\le 6m\tilde{n}\le n\quad (\omega\in\Omega).
\end{equation*}
and therefore
\begin{eqnarray}
\label{C2}
e_n^\ran(S^{N_1,N_2},B_{L_p^{N_1,N_2}},L_q^{N_1}) 
&\le & 
 c\Bigg(N_1^{1/p-1/q}\left\lceil\frac{\tilde{n}}{N_1}\right\rceil^{-\left(1-1/p\right)}+\left\lceil\frac{\tilde{n}}{N_1}\right\rceil^{-1/2}\Bigg).
\end{eqnarray}
Furthermore, using \eqref{C11}, we obtain
\begin{eqnarray}
\label{C3}
\left\lceil\frac{\tilde{n}}{N_1}\right\rceil 
&>&\frac{n}{12N_1\left\lceil c(1)\log(N_1+N_2)\right\rceil}
\ge \frac{n}{12N_1(c(1)+1)\log(N_1+N_2)}
\end{eqnarray}
and 
\begin{eqnarray*}
\left\lceil\frac{n}{N_1\log(N_1+N_2)}\right\rceil
&<&\frac{n+N_1\log(N_1+N_2)}{N_1\log(N_1+N_2)}\le \frac{n\left(1+\frac{1}{6c(1)}\right)}{N_1\log(N_1+N_2)}
\end{eqnarray*}
which together with \eqref{C3} yields
\begin{eqnarray}
\label{C5}
\left\lceil\frac{\tilde{n}}{N_1}\right\rceil 
&>&\frac{1}{12(c(1)+1)\left(1+\frac{1}{6c(1)}\right)}\left\lceil\frac{n}{N_1\log(N_1+N_2)}\right\rceil. 
\end{eqnarray}
Combining \eqref{C2} and \eqref{C5} completes the proof of  the upper bound in \eqref{A3}. Finally, the upper bound in the deterministic case \eqref{AC1} follows trivially from the norm bound \eqref{WN2}.

Now we prove the lower bounds in \eqref{A2}--\eqref{A4}. First assume $p\le 2$. Then the lower bound of  \eqref{A2} is a consequence of \eqref{UL2} and \eqref{UL3} of Proposition \ref{pro:1}. Next let $p>2$ and $p\ge q$. In this case the lower bound in \eqref{A2} follows from \eqref{UL1}. Now consider the case $2<p<q$. Here the lower bound of \eqref{A4} is a consequence of \eqref{UQ3}.  Finally, \eqref{UL1} and \eqref{UL2} imply  
\begin{eqnarray}
\lefteqn{  e_n^\ran(S^{N_1,N_2},B_{L_p^{N_1,N_2}},L_q^{N_1})  }
\notag\\
&\ge&c(2)N_1^{1/p-1/q}\left\lceil\frac{n}{N_1}\right\rceil^{-\left(1-1/p\right)}+c(2)\left\lceil\frac{n}{N_1}\right\rceil^{-1/2} \min \left( \log (N_1+1), \left\lceil\frac{n}{N_1}\right\rceil \right)^{\delta_{q,\infty}/2},\label{N9}
\end{eqnarray}
which in the case $q<\infty$ and in the case $(q=\infty) \wedge (\left\lceil n/N_1\right\rceil \ge \log (N_1+1))$ is just the lower bound in \eqref{A3}. Now assume $q=\infty$ and $\left\lceil n/N_1\right\rceil< \log (N_1+1)$. Then 
\begin{eqnarray*}
N_1^{1/p}\left\lceil\frac{n}{N_1}\right\rceil^{-\left(1-1/p\right)}&\ge& N_1^{1/p}(\log (N_1+1))^{-(1/2-1/p)}\left\lceil\frac{n}{N_1}\right\rceil^{-1/2}
\\
&\ge& c(3)\left\lceil\frac{n}{N_1}\right\rceil^{-1/2} (\log (N_1+1))^{1/2}.
\end{eqnarray*}
This combined with \eqref{N9} gives   
\begin{eqnarray*}
e_n^\ran(S^{N_1,N_2},B_{L_p^{N_1,N_2}},L_\infty^{N_1}) &\ge&c(2)N_1^{1/p}\left\lceil\frac{n}{N_1}\right\rceil^{-\left(1-1/p\right)}
\\
&\ge&\frac{c(2)}{2}N_1^{1/p}\left\lceil\frac{n}{N_1}\right\rceil^{-\left(1-1/p\right)}
+ \frac{c(2)c(3)}{2}\left\lceil\frac{n}{N_1}\right\rceil^{-1/2} (\log (N_1+1))^{1/2},
\end{eqnarray*}
thus the lower bound of \eqref{A3} also for that case. 

In the deterministic case we use Lemma \ref{lem:5} (iii) twice. We set $\bar{n}=N_1N_2$, $M=N_1$, consider the index-set $\Z[1,N_1]\times \Z[1,N_2]$, and let $n\in\N$, with  $4n<N_1N_2=\bar{n}$. First we treat the case $p\le q$.  Here we set $\psi_{ij}=N_1^{1/p}e_{ij}$, where $e_{ij}$ are the unit vectors in $\K^{N_1}\times \K^{N_2}$, and $\bar{I}_k =\{(k,j):\, 1\le j\le N_2\}$, thus $\bar{n}_k=N_2$ $(1\le k\le N_1)$. Furthermore, let
$$
n_k\in\N_0 \text{ with }n_k\le N_2 \text{ and }\sum_{k=1}^{N_1} n_k\le n,\quad I_k\subseteq \bar{I}_k, \text{ with }|I_k|= N_2-n_k.
$$
Since $n<\bar{n}/4$, it follows that there is a $k_0$ such that $n_{k_0}< N_2/4$, hence $|I_{k_0}|>3N_2/4$, which gives  
$$
\bigg\|\sum_{j:\, (k_0,j)\in I_{k_0}}S^{N_1,N_2}\psi_{k_0j}\bigg\|_{L_q^{N_1}}>\frac34 N_1^{1/p-1/q}.
$$
Now Lemma \ref{lem:5} (iii) together with \eqref{J9} implies \eqref{AC1} for $p\le q$. 

Next we assume $p>q$. Here we set $M=1$, $\psi_{ij}=e_{ij}$, $\bar{I}_1=\Z[1,N_1]\times \Z[1,N_2]$,  hence    $\bar{n}_1=\bar{n}=N_1N_2$. Furthermore, let $n\in\N$ with  $n<N_1N_2/4$, let $n_1\le n$ and $I_1\subseteq \bar{I}_1$ with $|I_1|= N_1N_2-n$. Then
$$
\bigg\|\sum_{(i,j)\in I_1}S^{N_1,N_2}\psi_{ij}\bigg\|_{L_q^{N_1}}\ge \bigg\|\sum_{(i,j)\in I_1}S^{N_1,N_2}e_{ij}\bigg\|_{L_1^{N_1}}=\frac{I_1}{N_1N_2}>\frac34,
$$
which together with \eqref{J9} completes the proof of \eqref{AC1}. 

\end{proof}

Let us have a look at the widest resulting gap between non-adaptive and adaptive randomized minimal errors in the region $N_1\le n<c(0)N_1N_2$, with $0<c(0)<1$ standing for the constant $c_0$ from Theorem \ref{theo:1}. Consider for $2<p<q$, $n\in\N$ 
\begin{equation*}
\gamma(p,q,n)=\max_{N_1,N_2:\,N_1\le n<c(0)N_1N_2}\; \frac{e_n^\ranno(S^{N_1,N_2},B_{L_p^{N_1,N_2}},L_q^{N_1})}{e_n^\ran(S^{N_1,N_2},B_{L_p^{N_1,N_2}},L_q^{N_1})}.
\end{equation*}
\begin{corollary}
\label{cor:1}
Let $2<p<q\le \infty$. Then there are constants $c_1,c_2>0$ such that for all $n\in\N$
\begin{equation}
\label{N1}
c_1n^{\frac{\left(\frac{1}{2}-\frac{1}{p}\right)\left(\frac{1}{p}-\frac{1}{q}\right)}{\frac{1}{2}-\frac{1}{q}}}(\log (n+1))^{-\left(1-1/p\right)}\le \gamma(p,q,n)\le c_2n^{\frac{\left(\frac{1}{2}-\frac{1}{p}\right)\left(\frac{1}{p}-\frac{1}{q}\right)}{\frac{1}{2}-\frac{1}{q}}}.
\end{equation}
\end{corollary}
\begin{proof}
It is convenient to estimate 
\begin{equation*}
\gamma(p,q,n)^{-1}=\min_{N_1,N_2:\,N_1\le n<c(0)N_1N_2}\; \frac{e_n^\ran(S^{N_1,N_2},B_{L_p^{N_1,N_2}},L_q^{N_1})}{e_n^\ranno(S^{N_1,N_2},B_{L_p^{N_1,N_2}},L_q^{N_1})}.
\end{equation*}
It follows from \eqref{A3} and \eqref{A4} of Theorem \ref{theo:1} that there are constants $c_1,c_2>0$ such that 
\begin{eqnarray}
&& c_1 \min_{N_1:\,N_1\le n}\max\bigg(\Big(\frac{n}{N_1}\Big)^{1/p-1/2}, 
N_1^{1/q-1/p} \bigg)
\le\gamma(p,q,n)^{-1}
\notag\\
&\le&c_2\min_{N_1,N_2:\,N_1\le n<c(0)N_1N_2}\left(\left(\Big(\frac{n}{N_1}\Big)^{1/p-1/2}+ 
N_1^{1/q-1/p} \right)(\log(N_1+N_2))^{1-1/p}\right)
\label{M5}
\end{eqnarray}
(for simplicity we omitted some log factors).
With $x_0$ satisfying
\begin{equation}
\Big(\frac{n}{x_0}\Big)^{1/p-1/2}= 
x_0^{1/q-1/p},
\label{N8}
\end{equation}
we have 
\begin{equation}
x_0=n^{\frac{\frac{1}{2}-\frac{1}{p}}{\frac{1}{2}-\frac{1}{q}}}, \quad x_0\in[1,n].\label{N7}
\end{equation}
and 
\begin{equation*}
\min_{x\in[1,n]}\max\left(\left(\frac{n}{x}\right)^{1/p-1/2}, x^{1/q-1/p}\right)= x_0^{1/q-1/p} =n^{-\frac{\left(\frac{1}{2}-\frac{1}{p}\right)\left(\frac{1}{p}-\frac{1}{q}\right)}{\frac{1}{2}-\frac{1}{q}}}.
\end{equation*}
This together with the lower bound in \eqref{M5} implies 
\begin{equation*}
cn^{-\frac{\left(\frac{1}{2}-\frac{1}{p}\right)\left(\frac{1}{p}-\frac{1}{q}\right)}{\frac{1}{2}-\frac{1}{q}}}
\le \gamma(p,q,n)^{-1}
\end{equation*}
and hence the upper bound in \eqref{N1}. 

Next we set 
\begin{equation*}
N_1=\left\lceil x_0\right\rceil,\quad 
N_2=\left\lfloor \frac{n}{c(0)x_0}\right\rfloor+1
\end{equation*}
implying 
\begin{equation*}
N_1\le n,\quad x_0\le N_1<2x_0,\quad \frac{n}{c(0)N_1}\le \frac{n}{c(0)x_0}<N_2<\frac{2n}{c(0)x_0}\le \frac{2n}{c(0)}\,,
\end{equation*}
so the requirement $N_1\le n<c(0)N_1N_2$ is fulfilled and the upper bound of  \eqref{M5} together with \eqref{N8} and  \eqref{N7} gives 
\begin{eqnarray*}
 \gamma(p,q,n)^{-1}
&\le&c\left(\left(\frac{n}{2x_0}\right)^{1/p-1/2}+ 
x_0^{1/q-1/p} \right)(\log(2x_0+2c(0)^{-1}n)^{1-1/p}
\\
&\le&c n^{-\frac{\left(\frac{1}{2}-\frac{1}{p}\right)\left(\frac{1}{p}-\frac{1}{q}\right)}{\frac{1}{2}-\frac{1}{q}}}
(\log(n+1))^{1-1/p},
\end{eqnarray*}
which yields the lower bound of \eqref{N1}.

\end{proof}
Consider the exponent of the  gap between non-adaption and adaption, for which we have
\begin{eqnarray*}
\frac{\left(\frac{1}{2}-\frac{1}{p}\right)\left(\frac{1}{p}-\frac{1}{q}\right)}{\frac{1}{2}-\frac{1}{q}}
&\le&\frac{\left(\frac{1}{2}-\frac{1}{q}\right)^2}{4\left(\frac{1}{2}-\frac{1}{q}\right)}=\frac{1}{8}-\frac{1}{4q}
\le \frac{1}{8},
\end{eqnarray*}
with equality everywhere iff $p=4$, $q=\infty$. With this choice the following holds.
For any $c_1,c_2$ with $c(0)^{1/2}<c_1<c_2$   a gap of order $n^{1/8}$ (up to log's) is reached for $N_1(n),N_2(n)\in [c_1 n^{1/2},c_2n^{1/2}]$ \; $(n\in\N, n\ge c_2^2)$.

\end{document}